\documentclass[final,a4paper,leqno]{siamltex}
\usepackage[english]{babel}
\usepackage[utf8]{inputenc}
\usepackage[labelformat=simple,font={it}]{subcaption}   
\usepackage[font={small,it}]{caption}
\usepackage{amsmath,amssymb,amsfonts,graphicx,outlines}
\usepackage{myshortcuts,float}
\usepackage{xcolor,bm}
\pagecolor{white}
\usepackage{fancyhdr}
\usepackage{MnSymbol}
\usepackage{csvsimple,mleftright}
\usepackage{tikz}
\usepackage{pgfplots}
\usetikzlibrary{decorations.pathreplacing}

\usepackage{appendix}
\usepackage{comment}
\newcommand{\norm}[1]{\left\lVert#1\right\rVert}
\usepackage[linesnumbered,lined,boxed,ruled,vlined]{algorithm2e}
\let\oldnl\nl
\newcommand{\nonl}{\renewcommand{\nl}{\let\nl\oldnl}}

\usepackage[hidelinks]{hyperref}
\usepackage{cleveref}
\newtheorem{remark}[theorem]{Remark}

\newenvironment{customlegend}[1][]{%
    \begingroup
    \csname pgfplots@init@cleared@structures\endcsname
    \pgfplotsset{#1}%
}{%
    \csname pgfplots@createlegend\endcsname
    \endgroup
}%
\def\addlegendimage{\csname pgfplots@addlegendimage\endcsname}

\title{Chebyshev HOPGD with sparse grid sampling for parameterized linear systems}
\author{Siobhán Correnty\thanks{Department of Mathematics, Royal Institute of Technology (KTH), SeRC Swedish
e-Science Research Center, Lindstedtsvägen 25, Stockholm, Sweden
\texttt{correnty@kth.se}} \and Melina A. Freitag\thanks{Institute of Mathematics, University of Potsdam, Karl-Liebknecht-Str. 24-25, 14476
Potsdam, Germany \texttt{melina.freitag@uni-potsdam.de}} \and Kirk M. Soodhalter\thanks{School of Mathematics, Trinity College Dublin, The University of Dublin, College Green, Dublin 2, Ireland \texttt{ksoodha@maths.tcd.ie}}}

\begin{document}

\maketitle

\begin{abstract}
We consider approximating solutions to parameterized linear systems of the form $A(\mu_1,\mu_2) x(\mu_1,\mu_2) = b$, where
	$(\mu_1, \mu_2) \in \mathbb{R}^2$. Here the matrix $A(\mu_1,\mu_2) \in \mathbb{R}^{n \times n}$ is nonsingular, large, and
	sparse and depends nonlinearly on the parameters $\mu_1$ and $\mu_2$. Specifically, the system arises from a
	discretization of a partial differential equation and $x(\mu_1,\mu_2) \in \mathbb{R}^n$, $b \in \mathbb{R}^n$.  The
	treatment of linear systems with nonlinear dependence on a single parameter has been well-studied, and robust methods
	combining companion linearization, Krylov subspace iterative methods, and Chebyshev interpolation have enabled fast
	solution for multiple values of a single parameter at the cost of a single iteration.  

	For systems depending nonlinearly on multiple parameters, the computation becomes much more challenging. This work
	overcomes those additional challenges by combining robust strategies for treating a single parameter with tensor-base
	approaches; i.e., combining companion linearization, the Krylov subspace method preconditioned bi-conjugate gradient (BiCG),
	and a decomposition of a tensor matrix of precomputed solutions, called snapshots. This produces a reduced order
	model of $x(\mu_1,\mu_2)$, and this model can be evaluated inexpensively for many values of the parameters. An
	interpolation of the model is used to produce approximations on the entire parameter space. In addition  this method can be used to
	solve a parameter estimation problem.

	This approach allows us to achieve similar computational savings as for the one-parameter case; we can solve for
	many parameter pairs at the cost of many fewer applications of an efficient iterative method.  The technique is presented
	for dependence on two parameters, but the strategy can be extended to more parameters using the same theoretical approach.
	Numerical examples of a parameterized Helmholtz equation show the competitiveness of our approach. The simulations are
	reproducible, and the software is available online. 
\end{abstract}



\section{Introduction}
We are interested in approximating solutions to linear systems of the form
\begin{align} \label{our-problem}
	A(\mu_1,\mu_2) x(\mu_1,\mu_2) = b,
\end{align}
for many different values of the parameters $\mu_1, \mu_2$. Here $A(\mu_1,\mu_2) \in \mathbb{R}^{n \times n}$ is a large and
sparse nonsingular matrix with a nonlinear dependence on $\mu_1 \in [a_1,b_1] \subset \mathbb{R}$ and $\mu_2 \in [a_2,b_2] \subset
\mathbb{R}$ and $x(\mu_1,\mu_2) \in \mathbb{R}^n$, $b \in \mathbb{R}^n$. Further, the system matrix arises from a discretization
of a parameterized partial differential equation (PDE) and can be expressed as the sum of products of matrices and functions,
i.e., 
\begin{align} \label{form-of-A}
	A(\mu_1,\mu_2) = A_0 + A_1 f_1(\mu_1,\mu_2) + \ldots + A_{n_f} f_{n_f} (\mu_1,\mu_2), 
\end{align}
where $n_f \ll n$ and $f_1,\ldots,f_{n_f}$ are nonlinear scalar functions in the parameters $\mu_1$ and $\mu_2$. This is usually referred to as \emph{split form}. 

The treatment of large-scale linear systems depending on a single parameter has been explored widely in the literature.
Initially, this was for the case of linear dependence on the parameter (e.g.,
\cite{ShiftedFreund,FrommerGlassner98,SOODHALTER2014105}). Strategies for treating nonlinear dependence on the parameter
\cite{CorrentyEtAl,CorrentyEtAl2,JarlebringCorrenty1} were then built on top of the strategies for treating the case of linear
dependence.  What these works
all have in common is that they allow for the evaluation of the solution for multiple values of the single parameter
at the cost of a single iteration of a Krylov subspace method.  This greatly reduces the costs, since the dominant costs (i.e.,
matrix-vector products and dot products) are independent from the number of evaluations of the solution for different parameter
values. Our work represents the next step: how to derive an approach for systems depending on multiple parameters.

The goal of achieving similar savings when the system depends on more than one parameter, such as in \eqref{our-problem}, is
much more challenging.  We describe how to build on the work for single-parameter-dependent linear systems in order achieve this
goal for the multi-parameter-dependent case. We combine the companion linearization-based Krylov subspace iteration employed in
\cite{CorrentyEtAl,CorrentyEtAl2,JarlebringCorrenty1} for the single parameter case with a tensor decomposition to construct a
reduced order model. This smaller model can be evaluated inexpensively to approximate the solution to \eqref{our-problem} for many
values of the parameters. Additionally, the model can be used to solve a parameter estimation problem, i.e., to simultaneously
estimate $\mu_1$ and $\mu_2$ for a given solution vector where these parameters are not known.  We describe our strategy in
terms of dependence of~\eqref{our-problem} on two parameters, but the theory extends to more than two parameters.

More specifically, our proposed method is based on a decomposition of a tensor matrix of precomputed solutions to
\eqref{our-problem}, called snapshots. In this way, building the reduced order model can be divided into two main parts:
\begin{enumerate}
	\item Generate the snapshots for many pairs $(\mu_1,\mu_2)$ efficiently using the technique from \cite{CorrentyEtAl2}.
	\item Perform the tensor decomposition.
\end{enumerate}

\section{From one to more parameters}\label{sec:single-multiple}
The transition from treating a system depending on a single parameter to one depending on multiple parameters presents an
interesting challenge.  The work in \cite{CorrentyEtAl,JarlebringCorrenty1} is based on generating a Taylor expansion of the
system matrix with respect to a single parameter centered around some point.  Generalizing this directly to the
multiple-parameter case can at times not produce satisfactory results due to a tight radius of convergence.  Instead, we find that
building on the approach in \cite{CorrentyEtAl2} allows for the derivation of a more robust algorithm for treating
\eqref{our-problem}.  Generating a truncated Taylor expansion is equivalent to constructing a Hermite interpolation at a single
point for a certain number of derivatives at that point.  We can interpret the work in \cite{CorrentyEtAl2} as replacing Hermite
interpolation with Chebyshev interpolation on the interval of interest for the single parameter. This leads to a different
companion linearization of the same problem, and one can evaluate for many values of parameters on the interval for the cost of a
single preconditioned Krylov subspace iteration (in this case a BiCG-based iteration \cite{BiCG76,Lanczos52}).  Furthermore,
Chebyshev interpolation comes with robust error estimates allowing for precise control of the interpolation error.

\subsection{Reduced order model from snapshots}\label{sec:rom-snapshots}

The strategy we propose to treat \eqref{our-problem} is to construct a reduced order model for $x(\mu_1,\mu_2)$ for which we
generate \emph{snapshots} (samples) varying one parameter at a time.  This allows us to efficiently acquire the snapshots using the robust
single-parameter algorithm proposed in \cite{CorrentyEtAl2}.
There have been many prior works that rely on \emph{sampling}; see, e.g.,
\cite{PODcit,HOPGD2,SparseHOPGD,HOPGD,HOPGD3,HarborAgitation}. In particular, solutions corresponding to each $(\mu_1,\mu_2)$ in
the tensor matrix are obtained by solving the linear systems individually. Using Preconditioned Chebyshev BiCG \cite{CorrentyEtAl2} to obtain the snapshots allows for great
flexibility in the selection of the set of snapshots included in the tensor, as the approximations are generated as one variable
functions of $\mu_1$ or of $\mu_2$, i.e., \emph{with one parameter fixed}. These functions are cheap to evaluate for different
values of the parameter, and this method is described in Section~\ref{sec:ChebyshevBiCG}. 

The tensor matrix of precomputed snapshots used here has a particular structure, as we consider sampling on a structured, sparse grid in the
parameter space. This is a way to overcome the so-called \textit{curse of dimensionality}, since the number of snapshots to
generate and store grows exponentially with the dimension when the sampling is performed on a conventional full grid. Note,
this work deals only with tensors of order $3$ in order to fully study the fundamental behavior of this method. 

The decomposition is performed with a variant of the higher-order proper generalized decomposition\footnote{The method HOPGD \cite{HarborAgitation} decomposes a tensor of snapshots sampled on a full grid.} (HOPGD), as proposed in
\cite{SparseHOPGD}. The
approach here has been adapted for this particular setting with fewer snapshots and, as in the standard implementation, results in
a separated expression in $\mu_1$ and $\mu_2$. Once constructed, the decomposition is interpolated to approximate solutions
to~\eqref{our-problem} corresponding to all parameters $(\mu_1,\mu_2) \in [a_1,b_1] \times [a_2,b_2]$. 

In general, we cannot guarantee that the error in a tensor decomposition for a certain set of snapshots will reach a specified
accuracy level. Further, it is not known a priori which sets of snapshots will lead to a successful decomposition, even with
modest standards for the convergence\footnote{We consider decompositions where the relative error at the nodes
\eqref{rel-err-nodes-footnote} is below $10^{-3}$ to have achieved modest convergence.} \cite{KoldaBader}. Different selections
of snapshots can lead to models that are more efficient to compute and more robust overall, but the optimal set is not known
beforehand. Moreover, generating the snapshots is the \emph{dominating cost} of building the reduced order model. In instances
where the decomposition process fails to converge or yields an inadequate model, our proposed method provides an effective remedy.
It facilitates the generation of a new set of snapshots within the same parameter space $[a_1, b_1] \times [a_2, b_2]$ with {\em
little extra computation}. These snapshots can be used to construct a new reduced order model for the same parameter space in an
identical way. This is the main contribution developed in this work; we refer to Section \ref{sec:novel} for additional
elaboration.

More precisely, an alternating directions, greedy algorithm is used to perform the decomposition, where the cost of each step of
the method grows linearly with the number of unknowns $n$. The basis vectors for the decomposition are not known beforehand,
similar to the established proper generalized decomposition (PGD), often used to create reduced order models of PDEs that are
separated in the time and space variables. The method PGD has been used widely in the study of highly nonlinear PDEs
\cite{PGD1,PGDoverview2014,PGD2} and has been generalized to solve multiparametric problems; see, for instance,
\cite{AmmarEtAl1,AmmarEtAl2}. 

Evaluating the resulting reduced order model is in general much more efficient than solving the systems individually for each
choice of the parameters \cite{SparseHOPGD}. Reliable and accurate approximations are available for a variety of parameter
choices, and parameter estimation can be performed in a similarly competitive way. Details regarding the construction of the
reduced order model are found in Section~\ref{sec:SparseHOPGD}, and numerical simulations of a parameterized Helmholtz equation
and a parameterized advection-diffusion equation are presented in Sections~\ref{sec:simulations} and \ref{sec:advec-diff},
respectively. All experiments were carried out on a 2.3 GHz Dual-Core Intel Core i5 processor with 16 GB RAM, and the
corresponding software can be found online.\footnote{\texttt{https://github.com/siobhanie/ChebyshevHOPGD}}

\subsection{Generating snapshots with Preconditioned Chebyshev BiCG} \label{sec:ChebyshevBiCG}
Reduced order models based on sampling of snapshot solutions are typically constructed by generating many finite element
approximations individually in an offline stage. We base our approach to obtain the snapshot solutions $x(\mu_1,\mu_2)$ to
\eqref{our-problem} on an adapted version of the Preconditioned Chebyshev BiCG method for parameterized linear systems, originally
proposed in \cite{CorrentyEtAl2}. Here, the unique solution to a companion linearization, formed from an arbitrarily accurate Chebyshev approximation
using \texttt{Chebfun} \cite{ChebfunGuide}, is generated in a preconditioned BiCG \cite{BiCG76,Lanczos52} setting.  Companion
linearization is common in the study of nonlinear eigenvalue problems. The linearization used here was chosen because of its prior
application to solving parameterized linear systems. Other choices may lead to competitive approaches, e.g.,
\cite{PEROVIC20181,DeTeran:2009:FIEDLER,VanBeeumen:2015:CORK}. 

Two executions of this method generate all finite element solutions corresponding to the values of the parameters shown in
Figure~\ref{tikz2}. Specifically, one execution for \emph{all solutions} on the line $\mu_1 = \mu_1^*$ and one execution for {\em
all solutions} on the line $\mu_2 = \mu_2^*$ in the plane $[a_1, b_1] \times [a_2, b_2]$. In this way, we generate solutions fixed
in one parameter. Moreover, the sampling represented here is sufficient to build a reliable reduced order model with the strategy
presented in \cite{SparseHOPGD}. We briefly summarize this method as follows.  

A Chebyshev approximation of $A(\mu_1,\mu_2)$ for a fixed $\mu_2^* \in \mathbb{R}$ is given by $P(\mu_1)$, where
\begin{align} \label{cheb1}
	 A(\mu_1,\mu_2^*) \approx P(\mu_1) \coloneqq P_0 \tau_0 (\mu_1) + \ldots + P_{d} \tau_{d} (\mu_1). 
\end{align}
Here $\tau_{\ell}$ is the degree $\ell$ Chebyshev polynomial on $[-c_1,c_1]$ and $P_{\ell}\in\mathbb{R}^{n \times n}$ are the corresponding interpolation coefficients. Note, $c_l \in \mathbb{R}$ is chosen such that $[a_l,b_l] \subseteq [-c_l,c_l]$, for $l=1,2$. A companion linearization based on the linear system
\begin{align} \label{cheby-approx}
	P (\mu_1) \hat{x} (\mu_1) = b
\end{align}
with $\hat{x}(\mu_1) = \hat{x}(\mu_1,\mu_2^*) \approx x(\mu_1,\mu_2^*)$, where $x(\mu_1,\mu_2^*)$ is as in \eqref{our-problem}, and a fixed $\mu_2^*$ has the form
\begin{align} \label{linear}
	\left( K - \mu_1 M \right) u(\mu_1) = c,
\end{align}
where $K, M \in \mathbb{R}^{d n \times d n}$ are coefficient matrices, independent of the parameter $\mu_1$, $c \in \mathbb{R}^{d n}$ is a constant vector, and the solution $u(\mu_1)$ is unique. This linearization, inspired by the work \cite{ChebyEffenKress} and studied fully in \cite{CorrentyEtAl2}, relies on the well-known three-term recurrence of the Chebyshev polynomials.

Solutions to the systems in \eqref{linear} are approximated with the Krylov subspace method bi-conjugate gradient for shifted systems \cite{PreCondAhmed,Frommer2003res}. Specifically, we approximate an equivalent shifted system, preconditioned with $(K-\sigma M)^{-1}$, given by
\begin{subequations} \label{precond}
\begin{eqnarray}
	&&( K - \mu_1 M )(K-\sigma M)^{-1} \hat{u}(\mu_1) = c \\
	&\iff& \left( I + (-\mu_1 + \sigma) M (K - \sigma M)^{-1} \right)\hat{u}(\mu_1) = c, \label{precond-2}
\end{eqnarray}
\end{subequations}
where $\hat{u}(\mu_1) = (K- \sigma M) u (\mu_1)$. The $k$th approximation comes from the Krylov subspace generated from the matrix $M (K - \sigma M)^{-1}$ and the vector $c$, defined by
\begin{align} \label{our-krylov}
	\mathcal{K}_k \coloneqq \text{span} \{ c, M (K - \sigma M)^{-1} c, \ldots, (M (K - \sigma M)^{-1} )^{k-1} c \},
\end{align}
where $\mathcal{K}_k = \mathcal{K}_k(M (K - \sigma M)^{-1},c)$. In particular, the method exploits the shift- and scaling- invariance properties of Krylov subspaces, and the preconditioner is applied efficiently based on the prior works \cite{Amiraslani2009,KressRoman}. 

A Lanczos biorthogonalization process generates a basis matrix $V_{k}\in\mathbb{R}^{d n \times k}$ for \eqref{our-krylov} and a tridiagonal matrix $T_k\in\mathbb{R}^{k \times k}$ such that 
\[
M(K - \sigma M)^{-1}V_k = V_k T_k + \beta_k v_{k+1} e_k^T 
\]
holds, where $\beta_k \in \mathbb{R}$, $I_k$ is the identity matrix of dimension $k$, and $e_k$ is the $k$th column of $I_k$. Further, $v_{k+1}$ is the $(k+1)$th column in the basis matrix for $\mathcal{K}_{k+1}$. 

Approximations to \eqref{linear} and, equivalently, approximations to \eqref{our-problem}, can be computed efficiently for $\mu_1^i$, $i=1,\ldots,n_1$, as follows:
\begin{subequations}\label{each-one-needs}
\begin{eqnarray}
	y_k (\mu_1^i) &=& \left(I_k + (-\mu_1^i + \sigma T_k)\right)^{-1} (\beta e_1),\label{solve-tri} \\
	z_k (\mu_1^i) &=& V_k y_k(\mu_1^i), \\
	x_k (\mu_1^i,\mu_2^*)  &=& \left( (K- \sigma M)^{-1} z_k(\mu_1^i) \right)_{1:n}.\label{krylov-approx}
\end{eqnarray}
\end{subequations}
Here, $\beta \coloneqq \norm{b}$, where $b$ is as in \eqref{our-problem}, and $x_k(\mu_1,\mu_2^*)$ denotes the approximation to \eqref{our-problem} from the Krylov subspace \eqref{our-krylov} of dimension $k$ corresponding to $(\mu_1,\mu_2^*)$. Equivalently, $x_k(\mu_1,\mu_2^*) \approx \hat{x}(\mu_1,\mu_2^*)$, where $\hat{x}(\mu_1,\mu_2^*)$ is the solution to the system \eqref{cheby-approx}. Note, the subscript on the right-hand side of \eqref{krylov-approx} denotes the first $n$ elements in the vector. 

Since $k \ll n$, solving the tridiagonal system in \eqref{solve-tri} is not computationally demanding. As a result, once we build a sufficiently large Krylov subspace via one execution of the main algorithm, we have access to approximations $x_k(\mu_1^i,\mu_2^*)$ for all $\mu_1^i$ in the interval $[a_1,b_1]$. The process of approximating \eqref{our-problem} for a fixed $\mu_1^*$ and $\mu_2^j$, for $j=1,\ldots,n_2$, is completely analogous to the above procedure and, thus, we have access to approximations $x_k(\mu_1^*,\mu_2^j)$ for all $\mu_2^j$ on the interval $[a_2,b_2]$, obtained in a similarly efficient way. 

\begin{remark}[Choice of target parameter] \label{remark-near-id}
The use of shift-and-invert preconditioners with well-chosen $\sigma$ generally result in fast convergence, i.e., a few iterations of the shifted BiCG algorithm. This is because $(K-\mu_1 M)(K-\sigma M)^{-1} \approx I$, for $\mu_1 \approx \sigma$. Note, taking $\sigma = \mu_1$ results in a breakdown of the method. 
\end{remark}

\section{Sparse grid based HOPGD} \label{sec:SparseHOPGD}
\begin{figure}[h]
\centering
\includegraphics[scale=1]{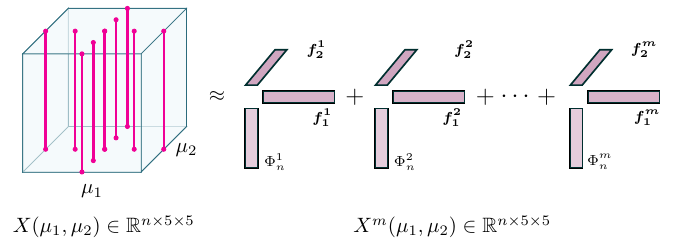}
\caption{Example of tensor matrix $X(\mu_1,\mu_2)$ consisting of 9 snapshot solutions (dots connected by vertical lines) and reduced order model $X^m(\mu_1,\mu_2)$ of rank $m$ as in \eqref{tensor}, where $\Phi_n^k \in \mathbb{R}^n$, $\bm{f_1^k}\coloneqq\begin{bmatrix} F_1^k(\mu_1^1),\ldots,F_1^k(\mu_1^5)\end{bmatrix} \in \mathbb{R}^5$, $\bm{f_2^k}\coloneqq\begin{bmatrix}F_2^k(\mu_2^1),\ldots,F_2^k(\mu_2^5) \end{bmatrix} \in \mathbb{R}^5$. Approximation is sum of rank-one tensors.}
\label{tikz1}
\end{figure}

\begin{figure}[h]
\centering
\includegraphics[scale=1]{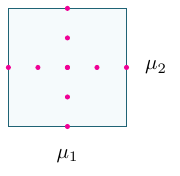}
\caption{Sampling on a sparse grid in the parameter space. Nodes \eqref{nodes} correspond to the 9 snapshot solutions in Figure~\ref{tikz1} with $\mu_1 \in [a_1,b_1]$, $\mu_2 \in [a_2,b_2]$, $\mu_1^* = (a_1+b_1)/2$, $\mu_2^* = (a_2 + b_2)/2$. All snapshots generated via two executions of Preconditioned Chebyshev BiCG.}
\label{tikz2}
\end{figure}

Let $X (\mu_1,\mu_2) \in \mathbb{R}^{n \times n_1 \times n_2}$ be a sparse, three-dimensional tensor of approximations to \eqref{our-problem}. These precomputed snapshots, generated by two executions of the preconditioned Krylov subspace method described in Section~\ref{sec:ChebyshevBiCG}, correspond to the pairs of parameters $(\mu_1^i,\mu_2^*)$, $i=1,\ldots,n_1$, and $(\mu_1^*,\mu_2^j)$, $j=1,\ldots,n_2$. Here, $\mu_1^i$ are in the interval $[a_1,b_1]$, $\mu_2^j$ are in the interval $[a_2,b_2]$, and $\mu_1^* \in [a_1,b_1]$, $\mu_2^* \in [a_2,b_2]$ are fixed values. More precisely, 
\begin{subequations}\label{snapshots}
\begin{alignat}{2} 
	X(\mu_1^i,\mu_2^*) &= \hat{x}(\mu_1^i,\mu_2^*) \in \mathbb{R}^n, \quad i&&=1,\ldots,n_1, \\
	X(\mu_1^*,\mu_2^j) &= \tilde{x}(\mu_1^*,\mu_2^j) \in \mathbb{R}^n,\quad j&&=1,\ldots,n_2,
\end{alignat}
\end{subequations}
and the remaining entries of $X(\mu_1,\mu_2)$ are zeros; see Figure~\ref{tikz1} for a visualization of a tensor of this form. Note, $\hat{x}(\mu_1^i,\mu_2^*)=\hat{x} (\mu_1^i)$ are approximations to the linear systems $A(\mu_1^i,\mu_2^*)x(\mu_1^i,\mu_2^*)=b$ and $\tilde{x}(\mu_1^*,\mu_2^j)=\tilde{x}(\mu_2^j)$ approximate solutions to $A(\mu_1^*,\mu_2^j)x(\mu_1^*,\mu_2^j)=b$. We refer to the set 
\begin{align} \label{nodes}
	\bm{\mu} \coloneqq \{ (\mu_1^i, \mu_2^*): i=1,\ldots,n_1 \} \cup \{ (\mu_1^*, \mu_2^j): j=1,\ldots,n_2 \}
\end{align}
as the \emph{nodes} and $\bm{\mu} \subset [a_1,b_1] \times [a_2,b_2]$. The set of nodes corresponding to the tensor matrix $X(\mu_1,\mu_2)$ in Figure~\ref{tikz1} are visualized in Figure~\ref{tikz2}. This way of sampling in the parameter space is a way to mitigate the so-called \textit{curse of dimensionality} in terms of the number of snapshots to generate and store.

Sparse grid based higher order proper generalized decomposition (HOPGD) \cite{SparseHOPGD} is a method which generates a tensor $X^m(\mu_1,\mu_2) \in \mathbb{R}^{n \times n_1 \times n_2}$ to approximate $X(\mu_1,\mu_2)$. Specifically, the approximations are separated in the parameters $\mu_1$ and $\mu_2$ and, at each pair of parameters in \eqref{nodes}, are of the form
\begin{subequations}\label{tensor}
\begin{alignat}{2} 
	X (\mu_1^i,\mu_2^*) \approx X^m (\mu_1^i,\mu_2^*) &= \sum_{k=1}^m \Phi_n^k F_1^k (\mu_1^i) F_2^k (\mu_2^*)  \in \mathbb{R}^{n},\quad i&&=1,\ldots,n_1, \\
	X (\mu_1^*,\mu_2^j) \approx X^m (\mu_1^*,\mu_2^j) &= \sum_{k=1}^m \Phi_n^k F_1^k (\mu_1^*) F_2^k (\mu_2^j) \in \mathbb{R}^{n},\quad j&&=1,\ldots,n_2,
\end{alignat} 
\end{subequations}
where $\Phi_n^k \in \mathbb{R}^n$ and 
\begin{align} \label{funcs}
	F_1^k: \mathbb{R} \rightarrow \mathbb{R}, \quad F_2^k: \mathbb{R} \rightarrow \mathbb{R}
\end{align}
are one variable scalar functions of the parameters $\mu_1$ and $\mu_2$, respectively. Here $m$ is the rank of the approximation and $F_{l}^k$, $l=1,2$, are the unknown functions of the $k$th mode. In this way, the nonzero entries of the tensor $X(\mu_1,\mu_2)$ are estimated by a linear combination involving products of lower-dimensional functions. 

The decomposition in \eqref{tensor} generates only the evaluations of $F_{l}^k$, and the reduced basis vectors for the approximation are not known a priori. Instead, these vectors are determined on-the-fly, in contrast to other methods based on the similar proper orthogonal decomposition \cite{PODcit,PGDoverview2014}. A visualization of this decomposition on a sparse tensor of snapshots appears in Figure~\ref{tikz1}. The spacing of the nodes is equidistant in both the $\mu_1$ and $\mu_2$ directions, conventional for methods of this type. 

We note the similarity between this method and the well-known higher-order singular value decomposition (HOSVD) \cite{HOSVDcite}. Additionally, the reduced order model generated via the method sparse grid based HOPGD is of the same form as the established CANDECOMP/PARAFAC (CP) decomposition. Specifically, the approximations consists of a sum of rank-one tensors \cite{KoldaBader}. This is visualized in Figure~\ref{tikz1}.

Efficient tensor decompositions have been the focus of several previous works. In \cite{doi:10.1137/060655894}, a Tucker-like approximation with was generated from an $n \times n \times n$ array in linear time, and a randomized algorithm was successfully implemented in order to compute a low-rank Tucker decomposition in \cite{Nakatsukasa2023}.
Moreover, in \cite{doi:10.1137/20M1382957}, a linear time version of a tensor adaptive cross approximation was presented. That method, which makes use of a rational approximation
algorithm, is fundamentally different from the method we propose to use, based on solving sequences of least squares approximations.

The model as expressed in \eqref{tensor} can only be evaluated on the set of nodes \eqref{nodes} corresponding to the snapshots included in the tensor matrix \eqref{snapshots}. To access other approximations, interpolations of $F_l^k$ will be required; see Section~\ref{sec:interpolation}. The particular structure of the sampling incorporated here is very similar to that which is used in traditional sparse grid methods for PDEs \cite{ActaSparse}. Our approach \emph{differs fundamentally} from these methods, which are based on an interpolation of high dimensional functions using a hierarchical basis. See \cite{HierBasisROM} for a method in this style, where a reduced order model was constructed using hierarchical collocation to solve parametric problems. 

\begin{remark} [Evaluation of the model at a point in space] 
The approximation in~\eqref{tensor} is expressed using vectors of length $n$ as the examples that follow in
	Sections~\ref{sec:simulations} concern a parameterized PDE discretized on a spatial domain denoted by $\Omega$. We are
	interested in the solutions to these PDEs for many values of the parameters $\mu_1$ and $\mu_2$. The approximation
	\eqref{tensor}, expressed at a point in space, is denoted by 
\begin{align}\label{one-point}
	X(x_0,\mu_1,\mu_2) \approx X^m (x_0,\mu_1,\mu_2) = \sum_{k=1}^m \Phi^k(x_0) F_1^k (\mu_1) F_2^k (\mu_2)  \in \mathbb{R},
\end{align}
for a particular $(\mu_1,\mu_2)\in\bm{\mu}$ \eqref{nodes}, where $x_0 \in \Omega \subset \mathbb{R}^d$ is a spatial parameter, $d$ is the dimension, and $\Phi^k(x_0) \in \mathbb{R}$ is an entry of the vector $\Phi_n^k \in\mathbb{R}^n$ in~\eqref{tensor}. Specifically, $\Phi_n^k$ denotes $\Phi^k$ evaluated on $\Omega_{dp}$, where $\Omega_{dp} \subset \Omega$. As in \eqref{tensor}, to evaluate the model in \eqref{one-point} for values of $\mu_1$ and $\mu_2$ outside of the set of nodes, interpolations must be performed.
\end{remark}

\subsection{Computation of the separated expression}
Traditionally, tensor decomposition methods require sampling performed on full grids in the parameter space, requiring the computation and storage of many snapshot solutions. See \cite{HOPGD2,HOPGD,HarborAgitation} for HOPGD performed in this setting, as well as an example in Section~\ref{sec:simulation-full}. The tensor decomposition proposed in \cite{SparseHOPGD} is specifically designed for sparse tensors. In particular, the sampling is of the same form considered in Section~\ref{sec:ChebyshevBiCG} and shown in Figures~\ref{tikz1} and \ref{tikz2}. This approach is more efficient, as it requires fewer snapshots and, as a result, the decomposition is performed in fewer computations. The method in \cite{SparseHOPGD}, adapted to our particular setting, is briefly summarized here.

We seek an approximation to the nonzero vertical columns of $X(\mu_1,\mu_2)$ described in \eqref{snapshots}, separated in the parameters $\mu_1 \in [a_1,b_1]$ and $\mu_2 \in [a_2,b_2]$. The full parameter domain is defined by $\mathcal{D} \coloneqq \Omega \times [a_1,b_1] \times [a_2,b_2]$, where $\Omega \subset \mathbb{R}^d$. We express solutions as vectors of length $n$, i.e., as approximations to \eqref{our-problem} evaluated at a particular node $(\mu_1,\mu_2)$; cf. \eqref{one-point}, where solutions are expressed at a single point in the spatial domain. The approximations satisfy
\begin{subequations}\label{approx-xm}
\begin{eqnarray}
	X(\mu_1,\mu_2) &\approx& X^m(\mu_1,\mu_2) \\ 
	 &=& \sum_{k=1}^m \Phi_n^k F_1^k (\mu_1) F_2^k (\mu_2) \label{second-line}\\
	&=& X^{m-1}(\mu_1,\mu_2) + \Phi_n^m F_1^m (\mu_1) F_2^m (\mu_2), \label{part-b}
\end{eqnarray}
\end{subequations}
for $(\mu_1,\mu_2) \in \bm{\mu}$, where $\bm{\mu}$ is the set of nodes defined in \eqref{nodes}, and are generated via an
alternating directions strategy. Specifically, at each enrichment step $m$, all components $\Phi_n^m \in \mathbb{R}^n$,
$F_1^m(\mu_1) \in \mathbb{R}$, and $F_2^m(\mu_2)\in \mathbb{R}$ in \eqref{part-b} are initialized and updated sequentially. The
components are used to construct an $L^2$ projection of $X(\mu_1,\mu_2)$, which is then added to the approximation. The updates
are done via a greedy algorithm, where all components are assumed fixed, except the one we seek to compute. The process is
equivalent to the minimization problem described as follows. At step $m$, find $X^m(\mu_1,\mu_2)  \in V_m \subset
L^2(\mathcal{D})$ as in \eqref{approx-xm} such that 
\begin{align}\label{min-prob}
	 \min_{X^m(\mu_1,\mu_2)} \left( \frac{1}{2} \norm{w X(\mu_1,\mu_2) - w X^m(\mu_1,\mu_2)}_{L^2 (\mathcal{D})}^2 \right), 
\end{align}
where $w$ is the sampling index, i.e.,
\[
\begin{cases}
w=1, & (\mu_1,\mu_2) \in \bm{\mu}, \\ w=0, & \text{otherwise},
\end{cases}
\]
$\mathcal{D}$ is the full parameter domain, and $V_m$ represents a set of test functions. Equivalently, using the weak formulation, we seek the solution $X^m(\mu_1,\mu_2)$ to 
\begin{align} \label{weak1}
	(wX^m(\mu_1,\mu_2), \delta X)_{\mathcal{D}} = (wX(\mu_1,\mu_2), \delta X)_{\mathcal{D}}, \quad \forall \delta X \in V_m.
\end{align}
Here $(\cdot,\cdot)_{\mathcal{D}}$ denotes the integral of the scalar product over the domain $\mathcal{D}$, and \eqref{weak1} can be written as
\begin{align} \label{solve}
	\left( w \Phi_n^m F_1^m F_2^m, \delta X \right)_{\mathcal{D}} = \left( w X(\mu_1,\mu_2) - wX^{m-1}(\mu_1,\mu_2), \delta X \right)_{\mathcal{D}}, \quad \forall \delta X \in V_m.
\end{align}
We have simplified the notation in \eqref{solve} with $F_1^m \coloneqq F_1^m(\mu_1)$ and $F_2^m \coloneqq F_2^m(\mu_2)$ for readability. In practice, approximations to \eqref{solve} are computed successively via a least squares procedure until a fixed point is detected. The $m$th test function is given by $\delta X \coloneqq \delta \Phi^m F_1^m F_2^m + \Phi_n^m \delta F_1^m F_2^m + \Phi_n^m F_1^m \delta F_2^m$, and the approximation \eqref{part-b} is updated with the resulting components. 

More concretely, let the initializations $\Phi_n^{m}$, $\bm{f_1^m} \coloneqq \begin{bmatrix} F_1^{m}(\mu_1^1),\ldots,F_1^m(\mu_1^{n_1}) \end{bmatrix}\in\mathbb{R}^{n_1}$, $\bm{f_2^m} \coloneqq \begin{bmatrix} F_2^{m}(\mu_2^1),\ldots,F_2^m(\mu_2^{n_2}) \end{bmatrix}\in\mathbb{R}^{n_2}$ be given. We seek first an update of $\Phi_n^{m_\text{old}}\coloneqq\Phi_n^{m}$, assuming $F_1^{m_{\text{old}}} \coloneqq F_1^{m}$, $F_2^{m_{\text{old}}} \coloneqq F_2^{m}$ are known function evaluations for all $(\mu_1,\mu_2) \in \bm{\mu}$. We update $\Phi_n^m$ as 
\begin{align}\label{phi-m}
	\Phi_n^m = \frac{\sum_{\bm{\mu}} r_{m-1}(\mu_1,\mu_2)} {\sum_{\bm{\mu}} F_1^m (\mu_1) F_2^m (\mu_2)},
\end{align}
where the $(m-1)$th residual vector is given by 
\begin{align} \label{res-u}
	r_{m-1}(\mu_1,\mu_2) \coloneqq X(\mu_1,\mu_2) - X^{m-1}(\mu_1,\mu_2) \in \mathbb{R}^n. 
\end{align}
Denote the vector obtained in \eqref{phi-m} by $\Phi_n^{m_{\text{new}}}$ and seek an update of $F_1^{m_\text{old}}(\mu_1^i)$, assuming $\Phi_n^m$, $F_2^m(\mu_2^j)$ known, using
\begin{align} \label{sol-f}
	F_1^m(\mu_1^i) = \frac{(\Phi_n^m)^T r_{m-1}(\mu_1^i,\mu_2^*)}{(\Phi_n^m)^T \Phi_n^m F_2^m(\mu_2^*)}, 
\end{align}
for $i=1,\ldots,n_1$. We denote the solutions found in \eqref{sol-f} by $F_1^{m_{\text{new}}}(\mu_1^i)$. Note, $\Phi_n^{m_\text{new}}$ and $F_2^{m_\text{old}}(\mu_2^j)$ are used in the computations in \eqref{sol-f}. The updates $F_2^{m}(\mu_2^j)$ are computed as
\begin{align} \label{sol-f2}
	F_2^m(\mu_2^j) = \frac{(\Phi_n^m)^T r_{m-1}(\mu_1^*,\mu_2^j)}{(\Phi_n^m)^T \Phi_n^m F_1^m(\mu_1^*)}, 
\end{align}
for $j=1,\ldots,n_2$, which we denote by $F_2^{m_\text{new}}(\mu_2^j)$, where $\Phi_n^{m_\text{new}}$ and $F_1^{m_{\text{new}}}(\mu_1^i)$ were used in the computations in \eqref{sol-f2}. For a certain tolerance level $\varepsilon_1$, if the relation 
\begin{align} \label{convergence1}
	\norm{\delta(\mu_1,\mu_2)} < \varepsilon_1
\end{align}
holds for all $(\mu_1,\mu_2) \in \bm{\mu}$, where 
\begin{align} \label{delta-mus}
	\delta (\mu_1,\mu_2) \coloneqq \Phi_n^{m_{\text{old}}} F_1^{m_{\text{old}}}(\mu_1) F_2^{m_{\text{old}}} (\mu_2)- \Phi_n^{m_{\text{new}}} F_1^{m_{\text{new}}}(\mu_1)  F_2^{m_{\text{new}}}(\mu_2),
\end{align}
we have approached a fixed point. In this case, the approximation in \eqref{part-b} is updated with the components $\Phi_n^m$, $F_1^m(\mu_1^i)$, and $F_2^m(\mu_2^j)$. If the condition \eqref{convergence1} is not met, we set $\Phi_n^{m_{\text{old}}} \coloneqq \Phi_n^{m_{\text{new}}}$, $F_1^{m_{\text{old}}} (\mu_1^i) \coloneqq F_1^{m_{\text{new}}}(\mu_1^i) $, $F_2^{m_{\text{old}}}(\mu_2^j) \coloneqq F_2^{m_{\text{new}}}(\mu_2^j)$ and repeat the process described above. If 
\begin{align} \label{rel-err-nodes-footnote}
	\frac{\norm{X(\mu_1,\mu_2) - X^{m}(\mu_1,\mu_2)}}{\norm{X(\mu_1,\mu_2) }} < \varepsilon_2
\end{align}
is met for a specified $\varepsilon_2$ and all $(\mu_1,\mu_2) \in \bm{\mu}$, the algorithm terminates. Otherwise, we seek $\Phi_n^{m+1}$, $F_1^{m+1}(\mu_1^i) $, $F_2^{m+1}(\mu_2^j) $ using the same procedure, initialized with the most recent updates. This process is summarized in Algorithm~\ref{alg:ChebyshevHOPGD}. 

In general, decompositions of many sets of snapshots can be performed to an acceptable accuracy level with $m \ll n$, where the parameter $m$ depends on the regularity of the exact solution \cite{PGDoverview2014}. Efficiency and robustness have been shown for the standard PGD method generated with a greedy algorithm like the one described here \cite{PGD2,HOPGD,Nouy10}. Note, the accuracy of the approximations here depends strongly on the quality of the separated functions in the model \cite{SparseHOPGD}. Additionally, the decomposition produced by the standard HOPGD method is optimal when the separation involves two parameters \cite{HarborAgitation}.  

The reduced order model described here can be used to approximate solutions to \eqref{our-problem} for many $(\mu_1, \mu_2)$
outside of the set of nodes. This is achieved through interpolating the one variable functions in \eqref{funcs}. Consider also a
solution vector $x (\mu_1,\mu_2)\in \mathbb{R}^n$ to \eqref{our-problem}, where $(\mu_1,\mu_2)$ are unknown. The interpolation of
$X^m(\mu_1,\mu_2)$ can be used to estimate the parameters $\mu_1$, $\mu_2$ that produce this particular solution. 

\begin{remark}[Separated expression in more than two parameters]
The reduced order model \eqref{tensor} is separated in the two parameters $\mu_1$ and $\mu_2$. In \cite{SparseHOPGD}, models separated in as many as six parameters were constructed, where the sampling was performed on a sparse grid in the parameter space. The authors note this particular approach for the decomposition may not be optimal in the sense of finding the best approximation, though it is possible. Our proposed way of computing the snapshots could be generalized to a setting separated in $s$ parameters. In this way, the procedure would fix $(s-1)$ parameters at a time and execute Chebyshev BiCG a total of $s$ times. 
\end{remark}

\begin{remark}[Comparison of HOPGD to similar methods]
In \cite{HarborAgitation}, HOPGD was studied alongside the similar HOSVD method. The decompositions were separated in as many as six parameters, and the sampling was performed on a full grid in the parameter space. In general, HOPGD produced separable representations with fewer terms compared with HOSVD. In this way, the model constructed by HOPGD can be evaluated in a more efficient manner. Furthermore, the method HOPGD does not require the number of terms in the separated solution to be set a priori, as in a CP decomposition. 
\end{remark}

\begin{algorithm}[h!]
\SetAlgoLined
\SetKwInput{Input}{Input}\SetKwInput{Output}{Output}\SetKwInput{Initialize}{Initialize}
\SetKwFor{FOR}{for}{do}{end}
\SetKwFor{IF}{if}{do}{end}
\SetKwFor{FUN}{function}{}{end}
\Input{Set of nodes $\bm{\mu}$ as in \eqref{nodes}, tensor matrix of precomputed snapshots $X(\mu_1,\mu_2)$ as in \eqref{tensor}, tolerances $\varepsilon_1,\varepsilon_2>0$}
\Output{Function $\bm{X}^m(\mu_1,\mu_2): [a_1,b_1]\times[a_2,b_2]\rightarrow\mathbb{R}^n$, approx to \eqref{our-problem}}
\Initialize{$\Phi_n^{1}\in\mathbb{R}^n$, $F_1^{1}(\mu_1^i)\in\mathbb{R}$, $F_2^{1}(\mu_2^j)\in\mathbb{R}$, $i=1,\ldots,n_1$, $j=1,\ldots,n_2$}
$m=1$ \\
Set $X^0(\mu_1^i,\mu_2^j) \coloneqq 0$, $X^m(\mu_1^i,\mu_2^j) \coloneqq \Phi_n^{m} F_1^{m}(\mu_1^i) F_2^{m}(\mu_2^j)$ \\
\While{$\norm{X(\mu_1,\mu_2) - X^{m}(\mu_1,\mu_2)}/\norm{X(\mu_1,\mu_2)} > \varepsilon_2$, for any $(\mu_1,\mu_2)\in\bm{\mu}$}{
Set $\Phi_n^{m_\text{old}} \coloneqq \Phi_n^{m}$, $F_1^{m_\text{old}}(\mu_1^i) \coloneqq F_1^{m}(\mu_1^i)$, $F_2^{m_\text{old}}(\mu_2^j) \coloneqq F_2^{m}(\mu_2^j)$ \\
Compute $\Phi_n^{m}$ with \eqref{phi-m} \\ 
Set $\Phi_n^{m_\text{new}}\coloneqq\Phi_n^{m}$\\
Compute $F_1^{m}(\mu_1^i)$ with \eqref{sol-f} \\ 
Set $F_1^{m_\text{new}}(\mu_1^i)\coloneqq F_1^{m}(\mu_1^i)$ \\
Compute $F_2^{m}(\mu_2^j)$ with \eqref{sol-f2} \\ 
Set $F_2^{m_\text{new}}(\mu_2^j)\coloneqq F_2^{m}(\mu_2^j)$ \\
Compute $\delta(\mu_1^i,\mu_2^j)$ using \eqref{delta-mus}  \\
\IF {$\delta(\mu_1,\mu_2) < \varepsilon_1$, for all $(\mu_1,\mu_2) \in \bm{\mu}$}
{
Update $X^m(\mu_1,\mu_2)$ as in \eqref{part-b} \\
Set $m \coloneqq m + 1$ \\ Set $\Phi_n^m \coloneqq \Phi_n^{m-1}$, $F_1^m (\mu_1^i) \coloneqq F_1^{m-1}(\mu_1^i)$, $F_2^m(\mu_2^j) \coloneqq F_2^{m-1}(\mu_2^j)$} 
}
{\textbf{return} Compute interpolation $\bm{X}^m(\mu_1,\mu_2)$ as in \eqref{int}} 
\caption{Chebyshev HOPGD with sparse grid sampling for parameterized linear systems}
\label{alg:ChebyshevHOPGD}
\end{algorithm}

\subsection{Interpolated model} \label{sec:interpolation}
The tensor representation in \eqref{tensor} consists, in part, of one-dimensional functions $F_1^k$ and $F_2^k$, for $k=1,\ldots,m$. Note, we have only evaluations of these functions at $\mu_1^i$, $i=1,\ldots,n_1$ and $\mu_2^j$, $j=1,\ldots,n_2$, respectively, and no information about these functions outside of these points. In this way, we cannot approximate the solution to \eqref{our-problem} for all $(\mu_1,\mu_2) \in [a_1,b_1] \times [a_2,b_2]$ using \eqref{tensor} as written. 

We can use the evaluations $F_1^k(\mu_1^i)$ and $F_2^k(\mu_2^j)$ in \eqref{tensor} to compute an interpolation of these one-dimensional functions in a cheap way. In practice, any interpolation is possible and varying the type of interpolation does not contribute significantly to the overall error in the approximation. Thus, we make the following interpolation of the representation in \eqref{tensor}:
\begin{align} \label{int}
	\bm{X}^m (\mu_1,\mu_2) &= \sum_{k=1}^m \Phi_n^k \bm{F}_1^k (\mu_1) \bm{F}_2^k (\mu_2)	
	 \in \mathbb{R}^{n},
\end{align} 
where $\bm{F}_1^k$, $\bm{F}_2^k$ are spline interpolations of $F_1^k$, $F_2^k$ \eqref{funcs}, respectively, and $\bm{X}^m(\mu_1,\mu_2): [a_1,b_1]\times[a_2,b_2]\rightarrow\mathbb{R}^n$. The interpolation in \eqref{int} can be evaluated to approximate~\eqref{our-problem} for other $(\mu_1, \mu_2) \in [a_1,b_1] \times [a_2,b_2]$ in a cheap way. This approximation to $x(\mu_1, \mu_2)$ is denoted as
\begin{align} \label{approx}
	x^m (\mu_1,\mu_2) \in \mathbb{R}^n,
\end{align}
and we compute the corresponding relative error as follows:
\[
	\frac{\norm{x^m(\mu_1,\mu_2) - x(\mu_1,\mu_2)}_2}{\norm{x(\mu_1,\mu_2)}_2} \cdot
\]
Note, simply interpolating several snapshots to estimate solutions to \eqref{our-problem} for all $(\mu_1,\mu_2)$ in the parameter space is not a suitable approach, as the solutions tend to be extremely sensitive to the parameters \cite{HarborAgitation}.

To solve the parameter estimation problem for a given solution $x(\mu_1,\mu_2) \in \mathbb{R}^n$ that depends on unknown parameters $(\mu_1,\mu_2)$, we use the \textsc{Matlab} routine \texttt{fmincon}. This routine uses the sequential quadratic programming algorithm \cite{NumerOptim} to find the pair of values $(\mu_1, \mu_2)$ in the parameter domain that minimizes the quantity
\begin{align} \label{min-prob2}
	\norm{x(\mu_1,\mu_2) - x^m (\mu_1,\mu_2)}_2. 
\end{align}
Simulations from a parameterized Helmholtz equation and a parameterized advection-diffusion equation appear in Sections~\ref{sec:simulations} and \ref{sec:advec-diff}, respectively. In these experiments, an interpolation of a reduced order model is constructed to approximate the solutions to the parameterized PDEs. 

\begin{remark}[Offline and online stages] \label{online-offline-remark}
In practice, the fixed point method in Algorithm~\ref{alg:ChebyshevHOPGD} can require many iterations until convergence, though the cost of each step scales linearly with the number of unknowns. Residual-based accelerators have been studied to reduce the number of these iterations \cite{SparseHOPGD}. This strategy, though outside the scope of this work, has shown to be beneficial in situations where the snapshots depended strongly on the parameters. Only one successful decomposition is required to construct \eqref{int}. Thus, the offline stage of the method consists of generating the snapshots with Chebyshev BiCG and executing Algorithm~\ref{alg:ChebyshevHOPGD}. 

Evaluating the reduced order model \eqref{int} requires only scalar and vector operations and, therefore, a variety of approximations to \eqref{our-problem} can be obtained very quickly, even when the number of unknowns is large. We consider this part of the proposed method the online stage. In the experiment corresponding to Figure~\ref{fig3}, evaluating the reduced order model had simulation time $0.010601$ CPU seconds. Generating the corresponding finite element solution with backslash in \textsc{Matlab} had simulation time $1.197930$ CPU seconds.
\end{remark}

\begin{remark}[Sources of error]
The approximations generated by Algorithm~\ref{alg:ChebyshevHOPGD} contain error from the Chebyshev approximation, the iterative method Preconditioned Chebyshev BiCG, the low-rank approximation, as well as the interpolations in \eqref{int}. The Chebyshev approximation, which is accurate to the level of machine precision, has a very small impact on the overall error, and Figure~\ref{fig1b} shows that the iterative method obtains relatively accurate approximations to the true solutions. In practice, the interpolations performed in~\eqref{int} are done on smooth functions, as visualized in Figure~\ref{fig-int}. The largest source of error from our proposed algorithm stems from the tensor decomposition, i.e., lines $1$--$17$ in Algorithm~\ref{alg:ChebyshevHOPGD}. This can be seen, for example, in Figure~\ref{fig2}. Conventional equidistant nodes construct the reduced order model from a variety of solutions. 
\end{remark}

\section{Novel contribution of our approach} \label{sec:novel}
Clearly, the most expensive operations in Algorithm~\ref{alg:ChebyshevHOPGD} are the inner products of two vectors of dimension
$n$, and, thus, the cost of each step in the tensor decomposition scales linearly with the number of unknowns in
\eqref{our-problem}. However, \emph{the choice of nodes \eqref{nodes} has a very large influence on the overall expense of the decomposition.} In particular, the sampling determines how many iterations are taken until a fixed point is found and how many modes $m$ are included in the reduced order model. If too many snapshots are included in the tensor matrix, the decomposition can fail to converge within expected time constraints. 

The choice of nodes can also determine the quality of the resulting model. Specifically, the decomposition considered in this work
is sensitive to \emph{overfitting}. When overfitting occurs, the snapshots are well approximated by the model, but approximations
corresponding to pairs of $(\mu_1,\mu_2)\notin\bm{\mu}$ are not acceptably accurate. Moreover, not including enough snapshots in
the tensor matrix can be problematic. This scenario can result in high levels of error, as the model has not received
sufficient information to capture the underlying structure of the solutions, a situation referred to as \emph{underfitting}.

It is not known a priori which sets of snapshots will produce acceptable models, as such convergence theory does not currently exist. Thus, it is advantageous to have access to solutions corresponding to many different pairs of $(\mu_1,\mu_2)$, as the initial selection rarely leads to a robust and accurate model. The traditional approach of generating snapshots one-by-one by solving finite element systems can be exceptionally costly. 

\emph{Our proposed method, which computes the snapshots as one-variable functions using Preconditioned Chebyshev BiCG, offers an effective means of generating a wide range of  accurate snapshots by utilizing a precomputed Krylov basis matrix. Consequently, accurate models for parameterized linear systems of the form \eqref{our-problem} can essentially be obtained by solving a tridiagonal system of dimension $k \times k$ for each snapshot.}

While the optimal set of snapshots is problem-dependent, our experience suggests that starting with a coarse selection of equally-spaced nodes and adding more until the overall error is reduced is a reasonable strategy. One must also take into account the storage demands related to the snapshots when the dimension $n$ is very large.  \color{black}

\section{Numerical simulations of a parameterized Helmholtz equation} \label{sec:simulations}
We consider both a reduced order model for a parameterized Helmholtz equation and a parameter estimation problem for the solution to a parameterized Helmholtz equation. Such settings occur naturally, for example, in the study of geophysics; see \cite{HarborAgitation,InverseHelmPar,GeophysicsCite}. These prior works were also based on a reduced order model, constructed using PGD. Similarly, the method PGD was used in the study of thermal process in \cite{ThermalProcessPGD}, and a reduced basis method for solving parameterized PDEs was considered in \cite{ReducedBasisMeth}. 

In the simulations which follow, the matrices $A_i$ arise from a finite element discretization, and $b$ is the corresponding load vector. All matrices and vectors here were generated using the finite element software FEniCS \cite{Alnaes:2015:FENICS}. The solutions to these discretized systems were approximated with a modified version of the Kryov subspace method Preconditioned Chebyshev BiCG \cite{CorrentyEtAl2}, as described in Section~\ref{sec:ChebyshevBiCG}. This strategy requires a linear solve with a matrix of dimension $n \times n$ on each application of the preconditioner and of its adjoint. We have chosen to perform the linear solve by computing one LU decomposition per execution of the main algorithm, which is reused at each subsequent iteration accordingly. This can be avoided by considering the inexact version of Preconditioned Chebyshev BiCG, described in the original work. 

This work proposes a novel improvement in generating snapshot solutions necessary for constructing models of this type. We choose to include three different examples in this section in order to fully capture the versatility of our strategy. Once the interpolation \eqref{int} has been constructed, evaluating the approximations for many values of the parameters $\mu_1$ and $\mu_2$ can be done in an efficient manner; see Remark~\ref{online-offline-remark}. In the case of the tensor decomposition failing to converge sufficiently or the model lacking accuracy outside of the original snapshots, a new set of approximations can be generated with little extra work. From these solutions, a new decomposition can be attempted; see Section~\ref{sec:novel}. 

\subsection{First simulation, snapshots sampled on a sparse grid} \label{sec:simulation-1}
Consider the Helmholtz equation given by
\begin{subequations} \label{model-prob}
\begin{alignat}{2}
\nabla^2 u(x) +f(\mu_1,\mu_2) u(x) &= h(x), \quad &&x \in \Omega, \\
u(x) &= 0, \quad &&x \in \partial \Omega,
\end{alignat}
\end{subequations}
where $\Omega \subset \mathbb{R}^2$ is as in Figure~\ref{fig3}, $\mu_1 \in [1,2]$, $\mu_2 \in [1,2]$, $f(\mu_1,\mu_2 ) = 2 \pi^2 + \cos(\mu_1) + \mu_1^4 + \sin(\mu_2) + \mu_2$, and $h(x) = \sin(\pi x_1) \sin(\pi x_2)$. A discretization of \eqref{model-prob} is of the form~\eqref{our-problem}, where
\[
	A(\mu_1,\mu_2) \coloneqq A_0 +  f(\mu_1,\mu_2) A_1. 
\]
We are interested in approximating $u(x)$ in \eqref{model-prob} for many different pairs of $(\mu_1,\mu_2)$ simultaneously. 

\begin{figure}[h]
	\centering
	\includegraphics[scale=.3]{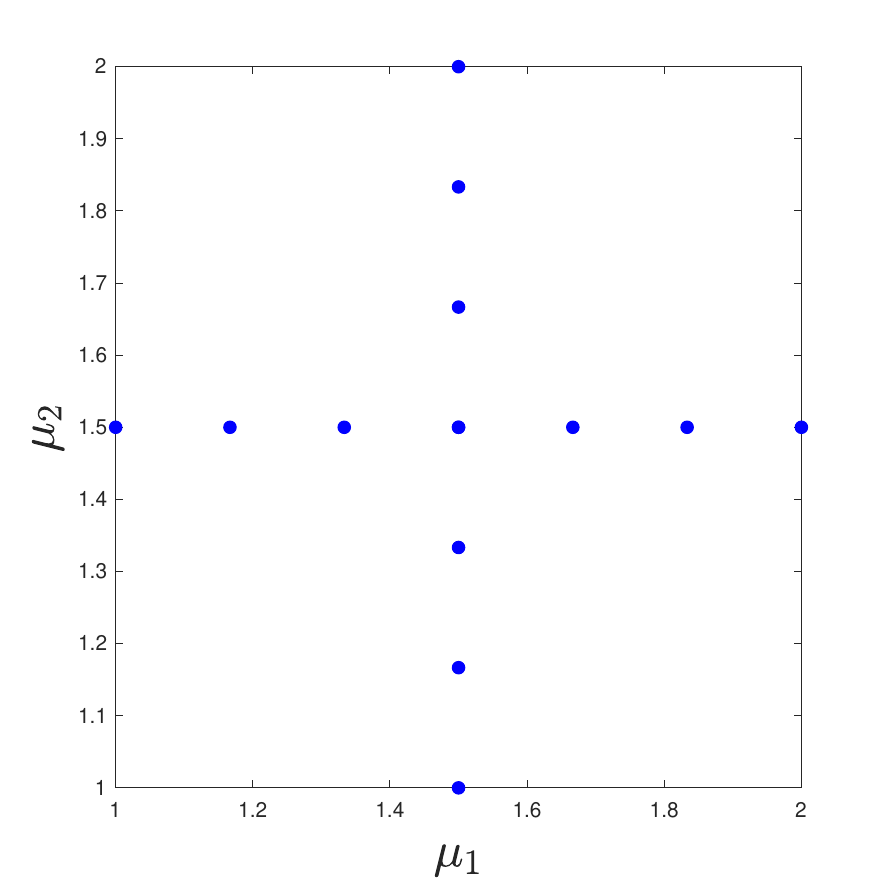}
	\caption{Locations of $13$ snapshots $(\mu_1^i,\mu_2^*)$ and $(\mu_1^*,\mu_2^j)$ used to construct $X^m$ in \eqref{tensor} to approximate \eqref{model-prob}, sampling on a sparse grid in the parameter space.}
	\label{fig1}
\end{figure}

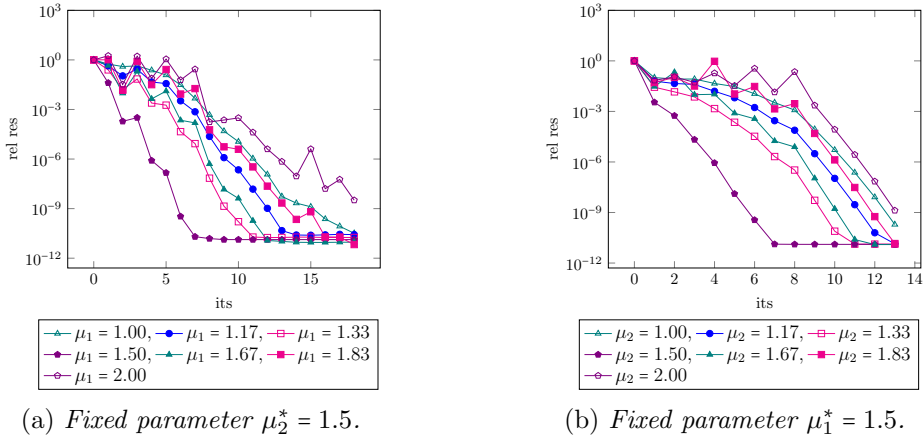
\begin{figure}[h]
	\centering

	\begin{subfigure}[T]{0.45\columnwidth}
		\centering
		\begin{tikzpicture}[scale=.6]
		\begin{semilogyaxis}[xlabel=its, ylabel=rel res,ymax=1.e3,legend style={at={(.45,-.2)},
		anchor=north,legend columns=3,font=\large},legend cell align={left},each nth point={1}]
	
		\addplot [mark=triangle,color=teal]
        		table[x=column 1,y=column 2,col sep=comma]{raw-data/sol1.csv};
		\addlegendentry{$\mu_1=1.00$,}
	
		\addplot [mark=*,color=blue]
        		table[x=column 1,y=column 2,col sep=comma]{raw-data/sol2.csv};
		\addlegendentry{$\mu_1=1.17$,}

		\addplot [mark=square,color=magenta]
        		table[x=column 1,y=column 2,col sep=comma]{raw-data/sol3.csv};
		\addlegendentry{$\mu_1=1.33$}
	
		\addplot [mark=pentagon*,color=violet]
        		table[x=column 1,y=column 2,col sep=comma]{raw-data/sol4.csv};
		\addlegendentry{$\mu_1=1.50$,}
		
		\addplot [mark=triangle*,color=teal]
        		table[x=column 1,y=column 2,col sep=comma]{raw-data/sol5.csv};
		\addlegendentry{$\mu_1=1.67$,}
		
		\addplot [mark=square*,color=magenta]
        		table[x=column 1,y=column 2,col sep=comma]{raw-data/sol6.csv};
		\addlegendentry{$\mu_1=1.83$}
		
		\addplot [mark=pentagon,color=violet]
        		table[x=column 1,y=column 2,col sep=comma]{raw-data/sol7.csv};
		\addlegendentry{$\mu_1=2.00$}

		\end{semilogyaxis}
		\end{tikzpicture}
	\caption{\centering Fixed parameter $\mu_2^* = 1.5$.}
	\label{chebyplot1}
	\end{subfigure}	
	\hfill
	\begin{subfigure}[T]{0.45\columnwidth}
		\centering
		\begin{tikzpicture}[scale=.6]
		\begin{semilogyaxis}[xlabel=its, ylabel=rel res,ymax=1.e3,legend style={at={(.45,-.2)},
		anchor=north,legend columns=3,font=\large},legend cell align={left},each nth point={1}]
	
		\addplot [mark=triangle,color=teal]
        		table[x=column 1,y=column 2,col sep=comma]{raw-data/sol1b.csv};
		\addlegendentry{$\mu_2=1.00$,}
	
		\addplot [mark=*,color=blue]
        		table[x=column 1,y=column 2,col sep=comma]{raw-data/sol2b.csv};
		\addlegendentry{$\mu_2=1.17$,}

		\addplot [mark=square,color=magenta]
        		table[x=column 1,y=column 2,col sep=comma]{raw-data/sol3b.csv};
		\addlegendentry{$\mu_2=1.33$}
	
		\addplot [mark=pentagon*,color=violet]
        		table[x=column 1,y=column 2,col sep=comma]{raw-data/sol4b.csv};
		\addlegendentry{$\mu_2=1.50$,}
		
		\addplot [mark=triangle*,color=teal]
        		table[x=column 1,y=column 2,col sep=comma]{raw-data/sol5b.csv};
		\addlegendentry{$\mu_2=1.67$,}
		
		\addplot [mark=square*,color=magenta]
        		table[x=column 1,y=column 2,col sep=comma]{raw-data/sol6b.csv};
		\addlegendentry{$\mu_2=1.83$}
		
		\addplot [mark=pentagon,color=violet]
        		table[x=column 1,y=column 2,col sep=comma]{raw-data/sol7b.csv};
		\addlegendentry{$\mu_2=2.00$}

		\end{semilogyaxis}
		\end{tikzpicture}
	\caption{\centering Fixed parameter $\mu_1^*=1.5$.}
	\label{chebyplot2}
	\end{subfigure}
	\caption{Convergence of Preconditioned Chebyshev BiCG for generating the snapshot solutions corresponding to the nodes in Figure~\ref{fig1} with $\sigma=1.48$.  Simulation in Figure~\ref{chebyplot1} gives all solutions for $(\mu_1,\mu_2^*)$, $\mu_1 \in [1,2]$, simulation in Figure~\ref{chebyplot2} gives all solutions for $(\mu_1^*,\mu_2)$, $\mu_2 \in [1,2]$.}
	\label{fig1b}
\end{figure}

The reduced order model $\bm{X}^m(\mu_1,\mu_2)$ is constructed as described in Algorithm~\ref{alg:ChebyshevHOPGD}. We can evaluate this model for $(\mu_1,\mu_2)$ outside of the nodes corresponding to the snapshot solutions. The particular $13$ nodes considered in this simulation are plotted in the parameter space, shown in Figure~\ref{fig1}. Figure~\ref{fig1b} displays the convergence of the two executions of Preconditioned Chebyshev BiCG required to generate the snapshots corresponding to the nodes in Figure~\ref{fig1}. Specifically, in Figure~\ref{chebyplot1}, we see faster convergence for approximations $\hat{x}(\mu_1,\mu_2^*)$ where the value $\mu_1$ is closer to the target parameter $\sigma$ in \eqref{precond}; see Remark~\ref{remark-near-id}. Note, we consider $\sigma=1.48$ to demonstrate the behavior of the method. An analogous result holds in Figure~\ref{chebyplot2}. We require the LU decompositions of $2$ different matrices of dimension $n \times n$ for this simulation. 

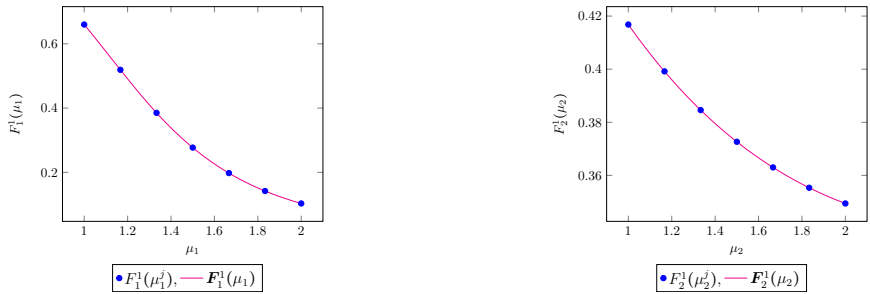
\begin{figure}[h]
	\centering
	\begin{subfigure}[T]{0.45\columnwidth}
		\centering
		\begin{tikzpicture}[scale=.5]
		\begin{axis}[xlabel=$\mu_1$, ylabel=$F_1^1(\mu_1)$,legend style={at={(.48,-.2)},
		anchor=north,legend columns=2,font=\large},legend cell align={left},each nth point={1}]
		
		\addplot [only marks, mark=*,color=blue]
        		table[x=column 1,y=column 2,col sep=comma]{raw-data/int2_f1.csv};
		\addlegendentry{$F_1^1(\mu_1^j)$, }
		
		\addplot [color=magenta]
        		table[x=column 1,y=column 2,col sep=comma]{raw-data/int_f1.csv};
		\addlegendentry{$\bm{F}_1^1(\mu_1)$}

		\end{axis}
		\end{tikzpicture}
	\end{subfigure}	
	\hfill
	\begin{subfigure}[T]{0.45\columnwidth}
		\centering
		\begin{tikzpicture}[scale=.5]
		\begin{axis}[xlabel=$\mu_2$, ylabel=$F_2^1(\mu_2)$,legend style={at={(.48,-.2)},
		anchor=north,legend columns=2,font=\large},legend cell align={left},each nth point={1}]
			
		\addplot [only marks, mark=*,color=blue]
        		table[x=column 1,y=column 2,col sep=comma]{raw-data/int2_f2.csv};
		\addlegendentry{$F_2^1(\mu_2^j)$, }
		
		\addplot [color=magenta]
        		table[x=column 1,y=column 2,col sep=comma]{raw-data/int_f2.csv};
		\addlegendentry{$\bm{F}_2^1(\mu_2)$}

		\end{axis}
		\end{tikzpicture}
	\end{subfigure}
	\caption{Function evaluations $F_1^1(\mu_1^i)$ and $F_2^1(\mu_2^j)$ produced by Algorithm~\ref{alg:ChebyshevHOPGD} to approximate~\eqref{model-prob}, plotted with corresponding interpolations $\bm{F}_1^1(\mu_1)$ and $\bm{F}_2^1(\mu_2)$ as in \eqref{int} with $\mu_1^i$ and $\mu_2^j$ as in Figure~\ref{fig1}.}
	\label{fig-int}
\end{figure}

\begin{figure}[h]
	\centering
	\begin{subfigure}[T]{0.45\columnwidth} \centering
	\includegraphics[scale=.32]{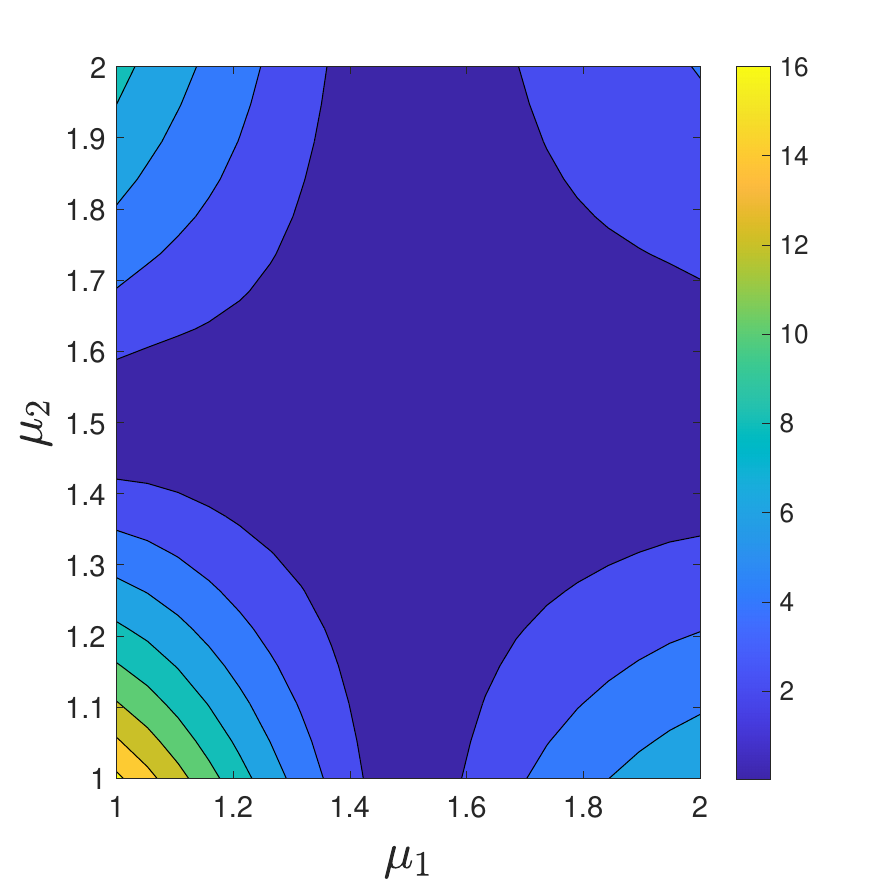}
	\label{}
	\end{subfigure}
	\hfill
	\begin{subfigure}[T]{0.45\columnwidth} \centering
	\includegraphics[scale=.32]{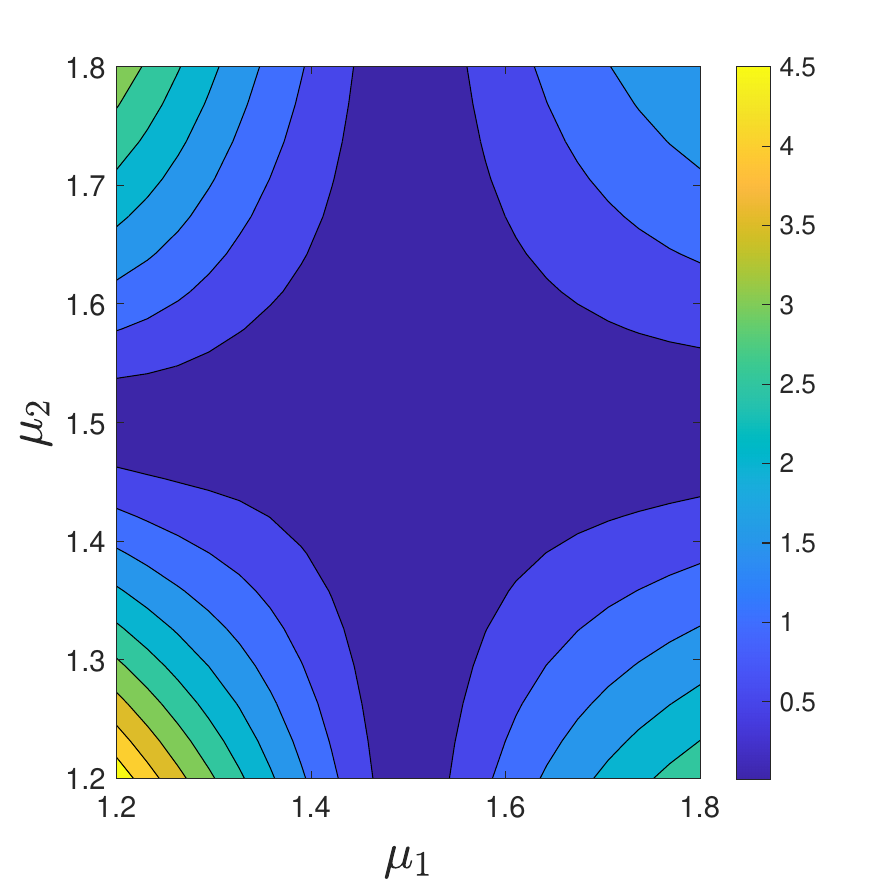}
	\label{}
	\end{subfigure}
	\caption{Percent relative error of $x^m(\mu_1,\mu_2)$ \eqref{approx} to approximate \eqref{model-prob} for 400 equidistant pairs of $(\mu_1,\mu_2)$, $m=14$, $n= 112995$. Here accurate: $< 1.5 \%$, reliable: $< 6 \%$. Note, same model in both figures.}
	\label{fig2}
\end{figure}

\begin{figure}[h]
\centering
\begin{subfigure}[T]{0.45\columnwidth}
\includegraphics[scale=0.4]{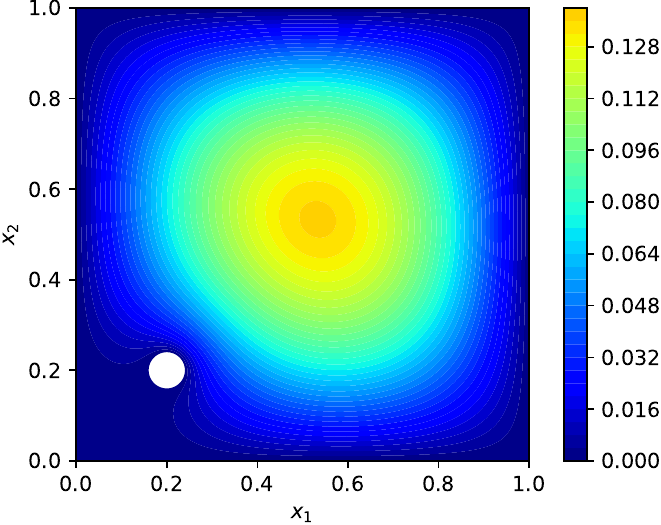}
\caption{Exact finite element solution}
\label{}
\end{subfigure}
\hfill
\begin{subfigure}[T]{0.45\columnwidth}
\includegraphics[scale=0.4]{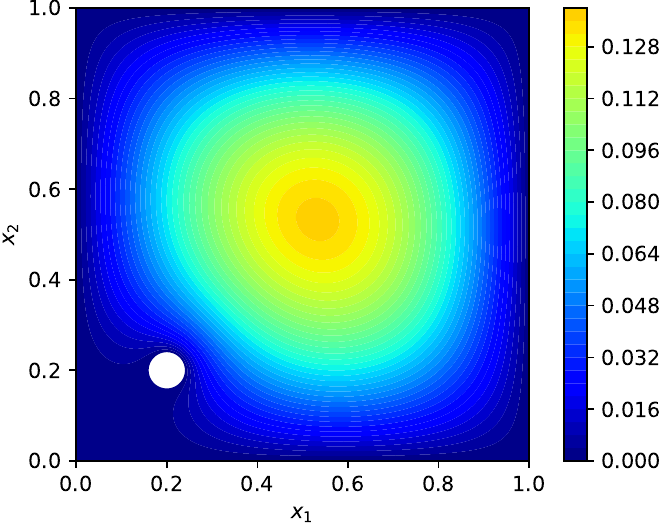}
\caption{Reduced order model solution}
\label{}
\end{subfigure}
\caption{Approximation to \eqref{model-prob} with $m=14$, $n=112995$, $0.347\%$ relative error (accurate), $\mu_1=1.6$, $\mu_2 = 1.4$.}
\label{fig3}
\end{figure}

\begin{figure}[h]
\centering
\begin{subfigure}[T]{0.45\columnwidth}
\includegraphics[scale=0.4]{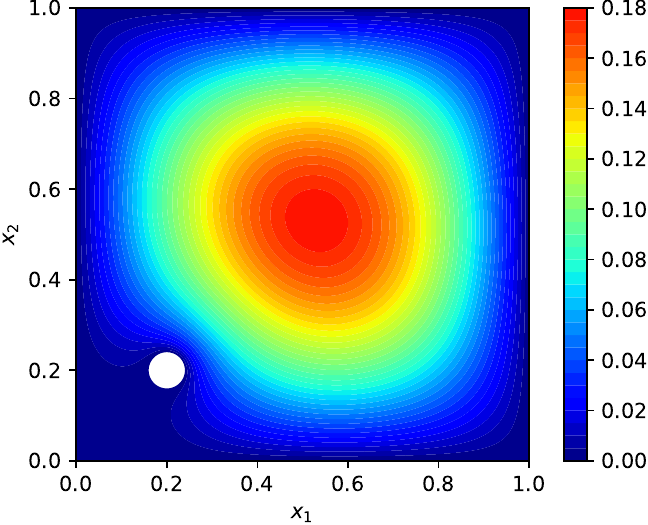}
\caption{Exact finite element solution}
\label{}
\end{subfigure}
\hfill
\begin{subfigure}[T]{0.45\columnwidth}
\includegraphics[scale=0.4]{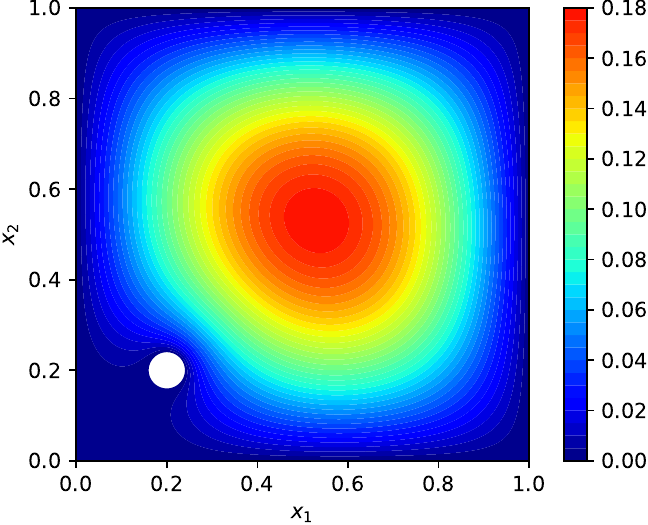}
\caption{Reduced order model solution}
\label{}
\end{subfigure}
\caption{Approximation to \eqref{model-prob} with $m=14$, $n=112995$, $0.171\%$ relative error (accurate), $\mu_1=1.45$, $\mu_2 = 1.6$.}
\label{fig4}
\end{figure}

\begin{figure}[h]
\centering
\begin{subfigure}[T]{0.45\columnwidth}
\includegraphics[scale=0.4]{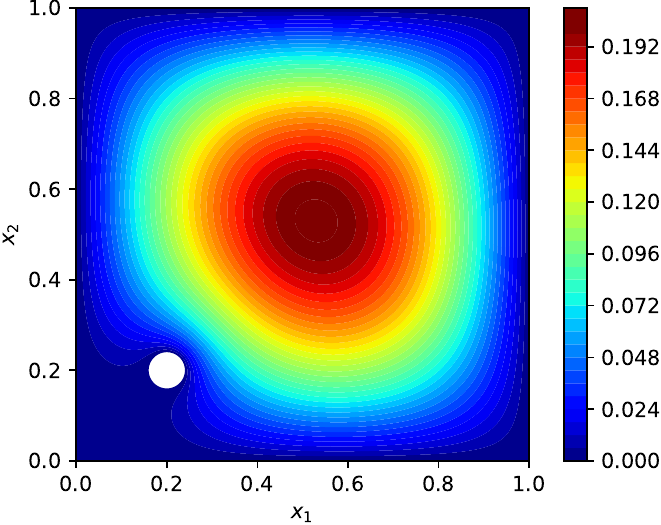}
\caption{Exact finite element solution}
\label{}
\end{subfigure}
\hfill
\begin{subfigure}[T]{0.45\columnwidth}
\includegraphics[scale=0.4]{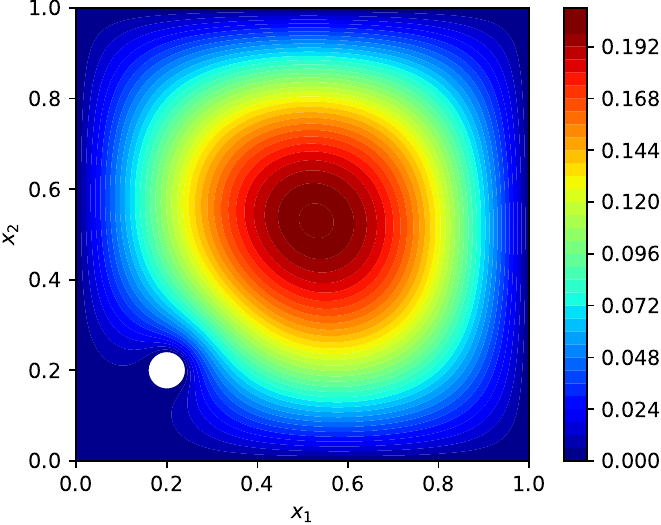}
\caption{Reduced order model solution}
\label{}
\end{subfigure}
\caption{Approximation to \eqref{model-prob} with $m=14$, $n=112995$, $0.414\%$ relative error (accurate), $\mu_1=1.4$, $\mu_2 = 1.4$.}
\label{fig5}
\end{figure}

In Figure~\ref{fig-int}, we see the interpolations $\bm{F}_1^1(\mu_1)$ and $\bm{F}_2^1(\mu_2)$ as in \eqref{int} plotted along with the function evaluations $F_1^1(\mu_1^i)$ and $F_2^1(\mu_2^j)$ generated by Algorithm~\ref{alg:ChebyshevHOPGD}, where $\mu_1^i$ and $\mu_2^j$ are as in Figure~\ref{fig1}. Figure~\ref{fig2} shows the percent relative error for approximations to \eqref{model-prob} for a variety of pairs $(\mu_1,\mu_2)$ different from the nodes plotted in Figure~\ref{fig1}. As in \cite{HOPGD}, we consider approximate solutions with percent relative error below $6 \%$ reliable and approximations with percent relative error below $1.5 \%$ accurate, where the error is computed by comparing the approximation to the solution obtained using backslash in \textsc{Matlab}. Figures~\ref{fig3}, \ref{fig4}, and \ref{fig5} show both the exact finite element solution and the solution generated from \eqref{int} for the sake of comparison. Here we display these solutions for $3$ pairs of $(\mu_1,\mu_2)$. These approximations all have percent relative error below $1.5 \%$. A variety of solutions can be obtained from this reduced order model. 

\subsection{Second simulation, snapshots sampled on a full grid} \label{sec:simulation-full}
To obtain an accurate model over the entire parameter space, we consider using HOPGD with sampling performed on a full grid. This is the approach taken in \cite{HOPGD2,HOPGD} to approximate parametric PDEs. The snapshots are obtained using the Krylov subspace method described in Section~\ref{sec:ChebyshevBiCG}, and the model is constructed in a way that is completely analogous to the method described in Section~\ref{sec:SparseHOPGD}. As the decomposition is performed on a tensor matrix containing more snapshots, the method is more costly. 

Specifically, we build a reduced order model to approximate the solution to the Helmholtz equation given by
\begin{subequations}\label{helmholtz1}
\begin{alignat}{2} 
	\left( \nabla^2 + f_1(\mu_1) + f_2(\mu_2) \alpha(x) \right) u(x) &= h(x), \quad &&x \in \Omega, \\
	u(x) &= 0, \quad &&x \in \partial \Omega,
\end{alignat}
\end{subequations}
where $\Omega \subset \mathbb{R}^2$ is as in Figure~\ref{fig7}, $\alpha(x) = x_2$, $\mu_1 \in [1,2]$, $\mu_2 \in [1,2]$, and $f_1(\mu_1) = \cos(\mu_1) + \mu_1^3$, $f_2(\mu_2) = \sin(\mu_2) + \mu_2^2$, $h(x)=\sin(\pi x_1) \sin(\pi x_2)$. A discretization of \eqref{helmholtz1} is of the form \eqref{our-problem}, where 
\[
	A(\mu_1,\mu_2) \coloneqq A_0 + f_1(\mu_1) A_1 + f_2 (\mu_2) A_2, 
\]
and approximating $u(x)$ in \eqref{helmholtz} for many different pairs of $(\mu_1,\mu_2)$ is of interest. We can still exploit the structure of the sampling when generating the snapshots by fixing one parameter and returning solutions as a one variable function of the other parameter. As in the simulations based on sampling on a sparse grid, if the decomposition fails to reach a certain accuracy level, or the model is inadequate due to overfitting or not enough data, a new set of snapshots can be generated with little extra computational effort. These solutions can be used to build a better model. 

\begin{figure}[h]
	\centering
	\begin{subfigure}[T]{0.45\columnwidth}
	\includegraphics[scale=.32]{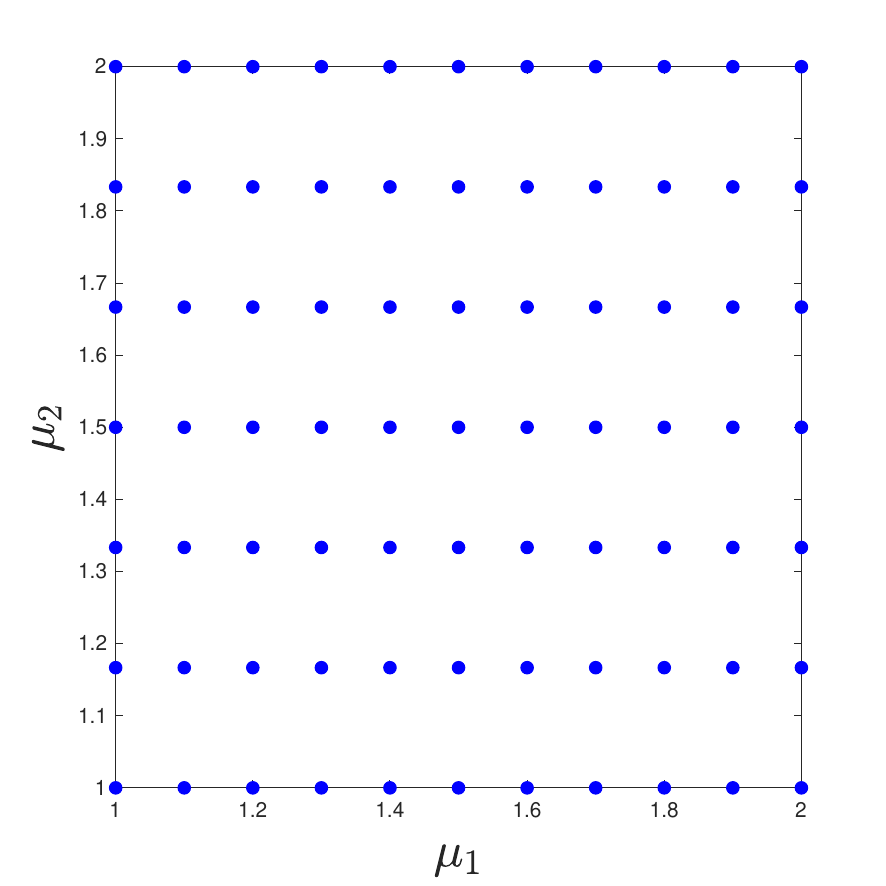}
	\caption{$77$ nodes, corresponding snapshots generated via $7$ executions of Preconditioned Chebyshev BiCG, $7$ fixed values $\mu_2^*$.}
	\label{fig6a}
	\end{subfigure}
	\hfill
	\begin{subfigure}[T]{0.45\columnwidth}
	\includegraphics[scale=.32]{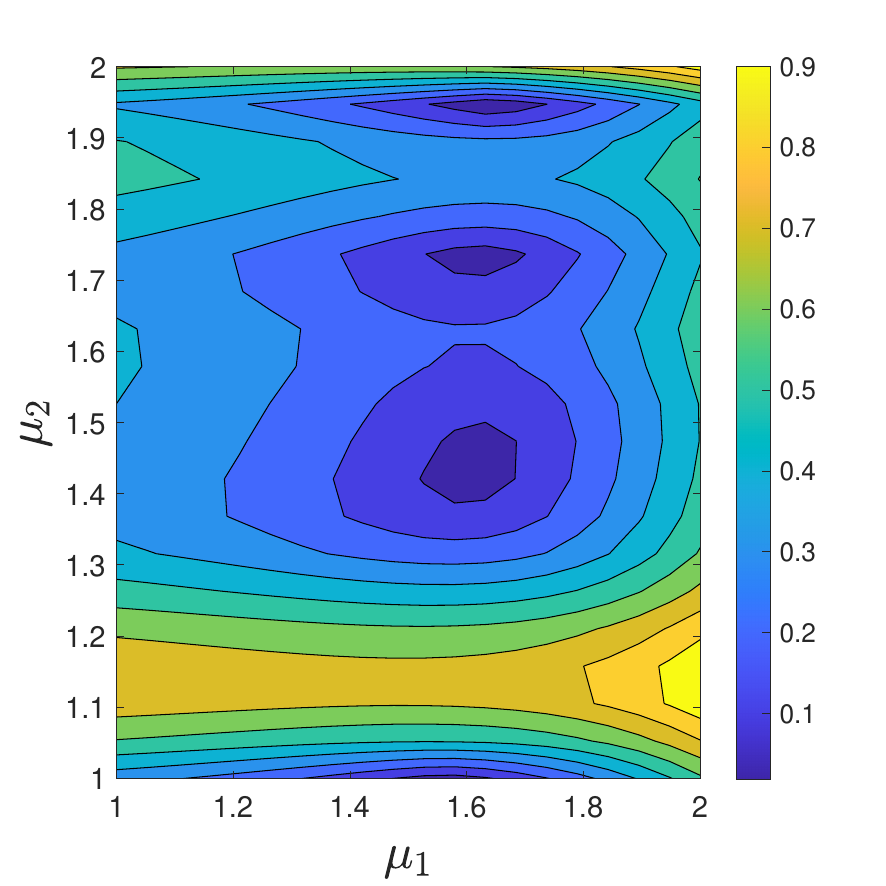}
	\caption{Percent relative error of interpolation of reduced order model for 400 equidistant pairs of $(\mu_1,\mu_2)$.}
	\label{fig6b}
	\end{subfigure}
	\caption{Simulation for approximating \eqref{helmholtz1}. Here percent relative error $< 1.5 \%$ (accurate), $n=112693$, $m=16$, sampling on a full grid.}
	\label{fig6}
\end{figure}

\begin{figure}[h]
	\centering
	\begin{subfigure}[T]{0.45\columnwidth}
	\includegraphics[scale=.4]{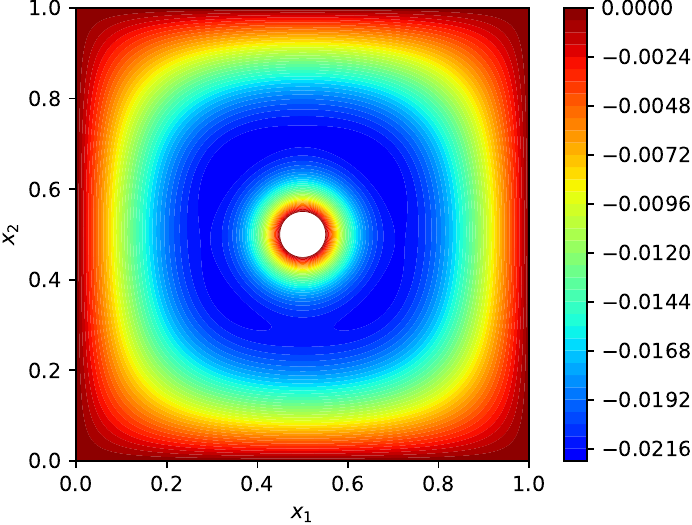}
	\label{}
	\caption{Exact finite element solution}
	\end{subfigure}
	\hfill
	\begin{subfigure}[T]{0.45\columnwidth}
	\includegraphics[scale=.4]{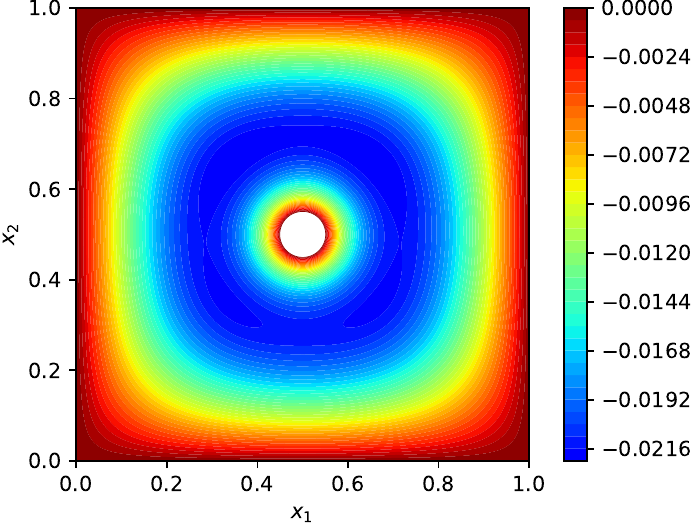}
	\label{}
	\caption{Reduced order model solution}
	\end{subfigure}
	\caption{Approximation to \eqref{helmholtz1} with $m=16$, $n=112693$, $0.340 \%$ relative error (accurate), $\mu_1=1$, $\mu_2=1$.}
	\label{fig7}
\end{figure}

\begin{figure}[h]
	\centering
	\begin{subfigure}[T]{0.45\columnwidth}
	\includegraphics[scale=.4]{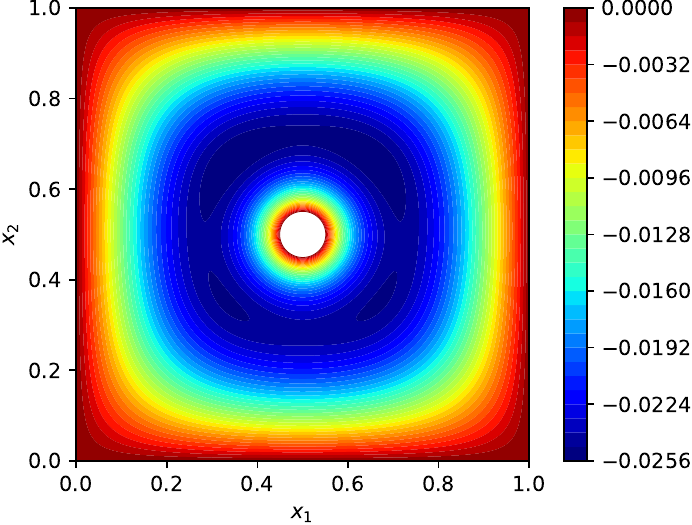}
	\label{}
	\caption{Exact finite element solution}
	\end{subfigure}
	\hfill
	\begin{subfigure}[T]{0.45\columnwidth}
	\includegraphics[scale=.4]{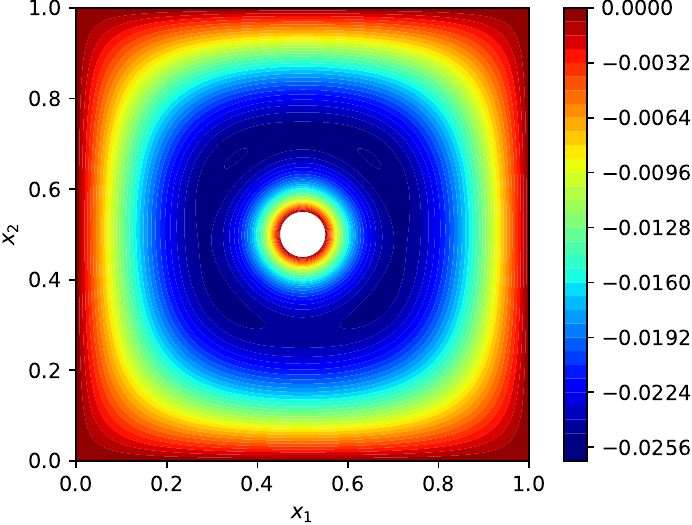}
	\label{}
	\caption{Reduced order model solution}
	\end{subfigure}
	\caption{Approximation to \eqref{helmholtz1} with $m=16$, $n=112693$, $0.714 \%$ relative error (accurate), $\mu_1=1.8$, $\mu_2=1.2$.}
	\label{fig8}
\end{figure}

\begin{figure}[h]
	\centering
	\begin{subfigure}[T]{0.45\columnwidth}
	\includegraphics[scale=.4]{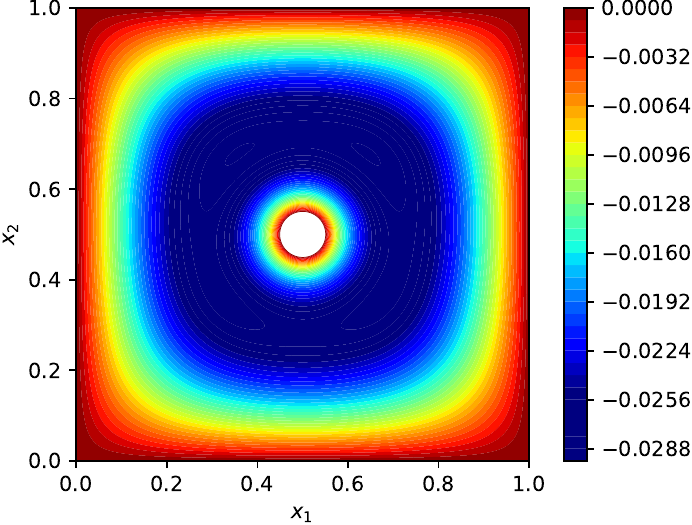}
	\label{}
	\caption{Exact finite element solution}
	\end{subfigure}
	\hfill
	\begin{subfigure}[T]{0.45\columnwidth}
	\includegraphics[scale=.4]{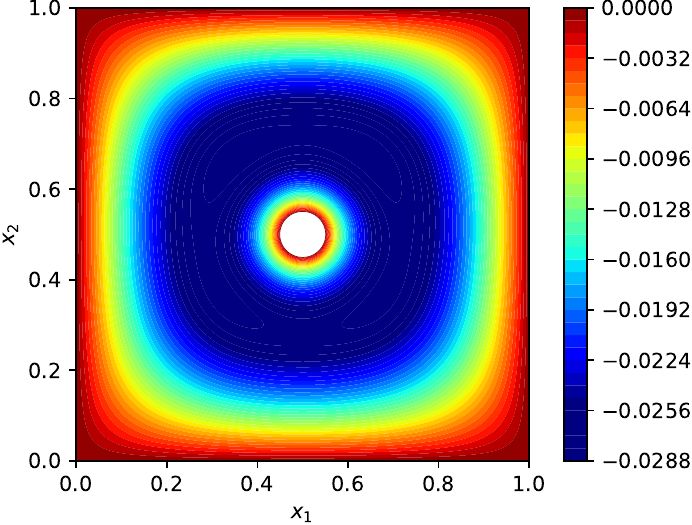}
	\label{}
	\caption{Reduced order model solution}
	\end{subfigure}
	\caption{Approximation to \eqref{helmholtz1} with $m=16$, $n=112693$, $0.946 \%$ relative error (accurate), $\mu_1=2$, $\mu_2=2$.}
	\label{fig9}
\end{figure}

Figure~\ref{fig6a} shows the locations of the nodes corresponding to $77$ different snapshot solutions used to construct a reduced order model. These snapshots were generated via $7$ executions of a modified form of Preconditioned Chebyshev BiCG, i.e., we consider $7$ different fixed $\mu_2^*$, equally spaced on $[a_2,b_2]$. This requires the LU decompositions of $7$ different matrices of size $n \times n$. Figure~\ref{fig6b} shows the percent relative error of the interpolated approximation, analogous to \eqref{int} constructed using Algorithm~\ref{alg:ChebyshevHOPGD}, for $400$ pairs of $(\mu_1,\mu_2) \in [1,2] \times [1,2]$, different from the nodes corresponding to the snapshot solutions. We note that the percent relative error of the approximation is below $1 \%$ for all pairs $(\mu_1,\mu_2)$, indicating that the approximations are accurate. For the sake of comparison, Figures~\ref{fig7}, \ref{fig8} and \ref{fig9} show both the exact finite element solutions and corresponding reduced order model solutions for $3$ pairs of $(\mu_1,\mu_2)$. In general, accurate solutions can be produced on a larger parameter space with sampling on a full grid.\footnote{See the adaptive refinement strategy proposed in \cite{SparseHOPGD}, used to achieve accuracy on the full parameter space.} Furthermore, the model can produce a variety of solutions. 

\subsection{Third simulation, parameter estimation with snapshots sampled on a sparse grid} \label{sec:param-est}
We consider now an application of our method to solve a parameter estimation problem. Similar methods have been constructed in the context of parameterized PDEs, for example, \cite{HOPGD3,ParEstSR,InverseHelmPar}. Our approach is analogous to the experiments performed in \cite{HOPGD3}, where several reduced order models were constructed using the method sparse grid based HOPGD \cite{SparseHOPGD}. Consider the Helmholtz equation given by
\begin{subequations}\label{helmholtz}
\begin{alignat}{2} 
	\left( \nabla^2 + f_1(\mu_1) + f_2(\mu_2) \alpha(x) \right) u(x) &= h(x), \quad &&x \in \Omega, \\
	u(x) &= 0, \quad &&x \in \partial \Omega,
\end{alignat}
\end{subequations}
where $\Omega = [0,1] \times [0,1]$, $\alpha(x) = x_1$, $f_1(\mu_1) = \sin^2 (\mu_1)$, $f_2(\mu_2) = \cos^2(\mu_2)$, $\mu_1 \in [0,1]$, $\mu_2 \in [0,1]$, and $h(x)=\exp(-x_1 x_2)$. A discretization of \eqref{helmholtz} is of the form \eqref{our-problem}, where 
\[
	A(\mu_1,\mu_2) \coloneqq A_0 + f_1(\mu_1) A_1 + f_2 (\mu_2) A_2.
\]
We consider a parameter estimation problem, i.e., we have a solution to the discretized problem, where the parameters $(\mu_1,\mu_2)$ are unknown. In particular, the solution is obtained using backslash in \textsc{Matlab}. 

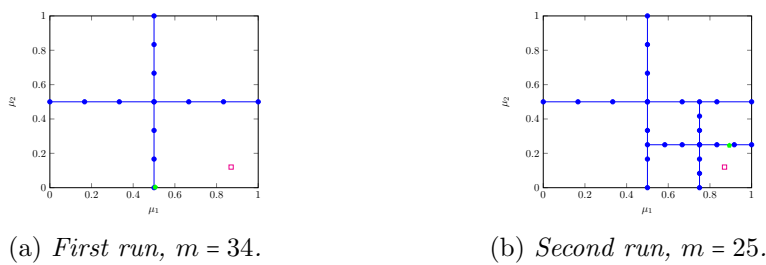
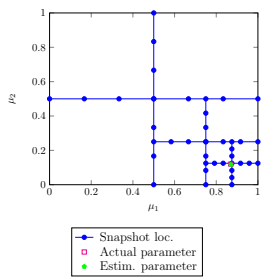
\begin{figure}[h]
\centering
		\begin{subfigure}[T]{0.45\columnwidth}
		\centering
		\begin{tikzpicture}[scale=.4]
		\begin{axis}[xlabel=$\mu_1$, ylabel=$\mu_2$,xmin=0,xmax=1,ymin=0,ymax=1,legend style={at={(.42,-.3)},
		anchor=north,legend columns=1,font=\large},legend cell align={left},each nth point={1}]

		\addplot [mark=*,color=blue]
        		table[x=column 1,y=column 2,col sep=comma]{raw-data/parest2.csv};

		\addplot [mark=pentagon*,color=green]
        		table[x=column 1,y=column 2,col sep=comma]{raw-data/parest4.csv};
		
		\addplot [mark=square,color=magenta]
        		table[x=column 1,y=column 2,col sep=comma]{raw-data/parest3.csv};
		
		\addplot [mark=*,color=blue]
        		table[x=column 1,y=column 2,col sep=comma]{raw-data/parest1.csv};

		\end{axis}
		\end{tikzpicture}
	\caption{\centering First run, $m=34$.}
	\label{fig4a}
	\end{subfigure}	
	\hspace{.4cm}
	\begin{subfigure}[T]{0.45\columnwidth}
		\centering
		\begin{tikzpicture}[scale=.4]
		\begin{axis}[xlabel=$\mu_1$, ylabel=$\mu_2$,xmin=0,xmax=1,ymin=0,ymax=1,legend style={at={(.45,-.2)},
		anchor=north,legend columns=1,font=\large},legend cell align={left},each nth point={1}]
	
		\addplot [mark=*,color=blue]
        		table[x=column 1,y=column 2,col sep=comma]{raw-data/parest6.csv};
	
		\addplot [mark=*,color=blue]
        		table[x=column 1,y=column 2,col sep=comma]{raw-data/parest1.csv};
	
		\addplot [mark=*,color=blue]
        		table[x=column 1,y=column 2,col sep=comma]{raw-data/parest2.csv};
		
		\addplot [mark=*,color=blue]
        		table[x=column 1,y=column 2,col sep=comma]{raw-data/parest5.csv};
	
		\addplot [mark=pentagon*,color=green]
        		table[x=column 1,y=column 2,col sep=comma]{raw-data/parest8.csv};
		
		\addplot [mark=square,color=magenta]
        		table[x=column 1,y=column 2,col sep=comma]{raw-data/parest7.csv};

		\end{axis}
		\end{tikzpicture}
	\caption{\centering Second run, $m=25$.}
	\label{fig4b}
	\end{subfigure}	
	
	\vspace{.3cm}
	
	\begin{subfigure}[T]{0.45\columnwidth}
		\centering
		\begin{tikzpicture}[scale=.4]
		\begin{axis}[xlabel=$\mu_1$, ylabel=$\mu_2$,xmin=0,xmax=1,ymin=0,ymax=1,legend style={at={(.42,-.25)},legend entries={abc},scatter/classes={
		abc={only marks}},anchor=north,legend columns=1,font=\large},legend cell align={left},each nth point={1}]
	
		\addplot [mark=*,color=blue]
        		table[x=column 1,y=column 2,col sep=comma]{raw-data/parest9.csv};
		\addlegendentry{\hphantom{ }Snapshot loc.}	
		
		\addplot [mark=square,color=magenta,only marks]
        		table[x=column 1,y=column 2,col sep=comma]{raw-data/parest11.csv};
		\addlegendentry{\hphantom{ }Actual parameter}
		
		\addplot [mark=pentagon*,color=green,only marks]
        		table[x=column 1,y=column 2,col sep=comma]{raw-data/parest12.csv};
		\addlegendentry{\hphantom{ }Estim. parameter}
	
		\addplot [mark=*,color=blue]
        		table[x=column 1,y=column 2,col sep=comma]{raw-data/parest1.csv};
	
		\addplot [mark=*,color=blue]
        		table[x=column 1,y=column 2,col sep=comma]{raw-data/parest2.csv};
		
		\addplot [mark=*,color=blue]
        		table[x=column 1,y=column 2,col sep=comma]{raw-data/parest5.csv};
	
		\addplot [mark=*,color=blue]
        		table[x=column 1,y=column 2,col sep=comma]{raw-data/parest6.csv};
		
		\addplot [mark=*,color=blue]
        		table[x=column 1,y=column 2,col sep=comma,skip coords between index={3}{4}]{raw-data/parest10.csv};
	
		\end{axis}
		\end{tikzpicture}
	\caption{\centering Third run, $m=7$.}
	\label{fig4c}
	\end{subfigure}	
\caption{Parameter estimation for solution to \eqref{helmholtz}, $n=10024$.}
\label{}
\end{figure}

Figures~\ref{fig4a}, \ref{fig4b}, and \ref{fig4c} together show a method similar to the one proposed in \cite{HOPGD3} for approximating this problem, performed by constructing $3$ successive HOPGD models with sampling on sparse grids. This simulation executes a modified form of Preconditioned Chebyshev BiCG $6$ times, requiring the LU factorizations of $6$ different matrices of dimension $n \times n$ to generate the $39$ snapshot solutions used to create the $3$ reduced order models. Once the models have been constructed, an interpolation is performed, followed by the minimization of \eqref{min-prob2}. Table~\ref{table:1} shows the percent relative error in the estimated values of $\mu_1$ and $\mu_2$ for each run. Each execution of our strategy leads to a better estimation of the pair of parameters. As before, if the decomposition fails to reach a certain level of accuracy on a set of snapshots, or the model is lacking in robustness, a new set of approximations in the same parameter space can be generated with little extra computational effort. These solutions can be used in order to construct a better model. 

\begin{table}[H]
\centering
\begin{center}
\begin{tabular}{ c | r | r }
 & rel err $\%$ $\mu_1$ & rel err $\%$ $\mu_2$ \\ \hline
First run, $m=34$ & 41.87417 & 97.40212 \\ \hline  
Second run, $m=25$ & 2.66555 & 105.09263 \\ \hline  
Third run, $m=7$ & 0.05706 & 1.69921 \\ \hline
\end{tabular}
\end{center}
\caption{Percent relative error for parameter estimation in $\mu_1$ and $\mu_2$ for solution to \eqref{helmholtz}, $n=10024$.}
\label{table:1}
\end{table}

\section{Numerical simulations of a parameterized advection-diffusion equation} \label{sec:advec-diff}
The advection-diffusion equation can be used to model particle transport, for example the distribution of air pollutants in the atmosphere; see \cite{1d11d082-9fb9-3fdc-a755-55dde241b26e,Ulfah_2018}. In particular, we consider the parameterized advection-diffusion equation given by
\begin{align} \label{ad-diff}
	\frac{\partial}{\partial t} u(x,t) = f_1(\mu_1) \frac{\partial^2}{\partial x^2} u(x,t) + f_2(\mu_2) \frac{\partial}{\partial x}u(x,t),
\end{align}
where $f_1(\mu_1) = 1 + \sin(\mu_1)$, $f_2 (\mu_2) = 10 + \cos(\mu_2) + \pi \mu_2$ and $x\in[0,1]$, $t \in [0,T]$. The boundary conditions $u(0,t)=0$, $u(1,t)=0$ are enforced, as well as the initial condition $u_0(x) = u(x,0)=\sin (\pi x)$. Here, $f_1(\mu_1)$ is referred to as the diffusion coefficient, and the value $f_2(\mu_2)$ is the advection parameter. Discretizing~\eqref{ad-diff} in space, with a finite difference scheme, and in time, with the implicit Euler method, gives a parameterized linear system of the form \eqref{our-problem}, where
\[
	A(\mu_1,\mu_2) \coloneqq A_0 + f_1(\mu_1) A_1 + f_2(\mu_2) A_2. 
\]
Note, the solution to this linear system gives an approximation to \eqref{ad-diff} at a specific discrete time-step. As in Section~\ref{sec:simulations}, we construct a reduced order model via a tensor decomposition \cite{SparseHOPGD}, where the snapshots are generated efficiently with a modified version of the method Preconditioned Chebyshev BiCG \cite{CorrentyEtAl2}. Here the sampling is performed on a sparse grid in the parameter space. If the tensor decomposition is not successful or the resulting model is not sufficiently accurate, a new set of snapshots can be generated on the same parameter space with little extra computation. These solutions can be used to attempt to build a better model. Similar model problems appear in, for example, \cite{Nouy10}, where the method PGD was used to construct approximate solutions to the parameterized PDE. Once the corresponding reduced order model \eqref{int} has been constructed, approximating the solutions for many $(\mu_1,\mu_2)$ can be done very cheaply; see Remark~\ref{online-offline-remark}. 

\subsection{First simulation, snapshots sampled on a sparse grid}
\begin{figure}[h]
	\centering
	\begin{subfigure}[T]{0.45\columnwidth}
	\includegraphics[scale=.32]{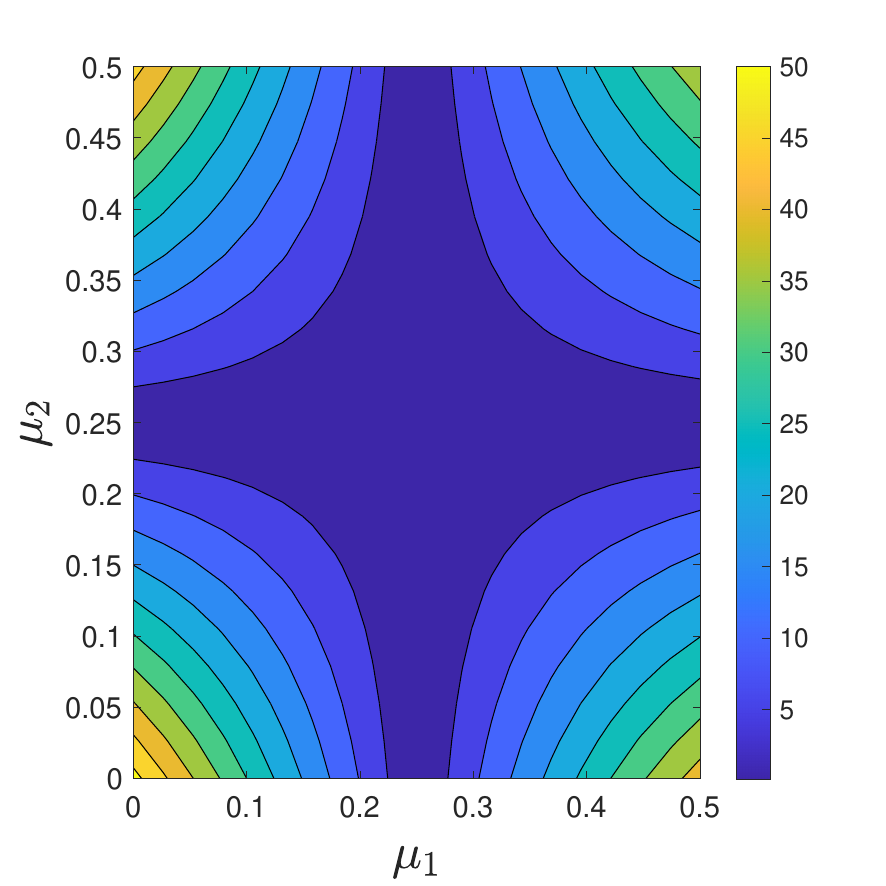}
	\label{fig10a}
	\end{subfigure}
	\hfill
	\begin{subfigure}[T]{0.45\columnwidth}
	\includegraphics[scale=.32]{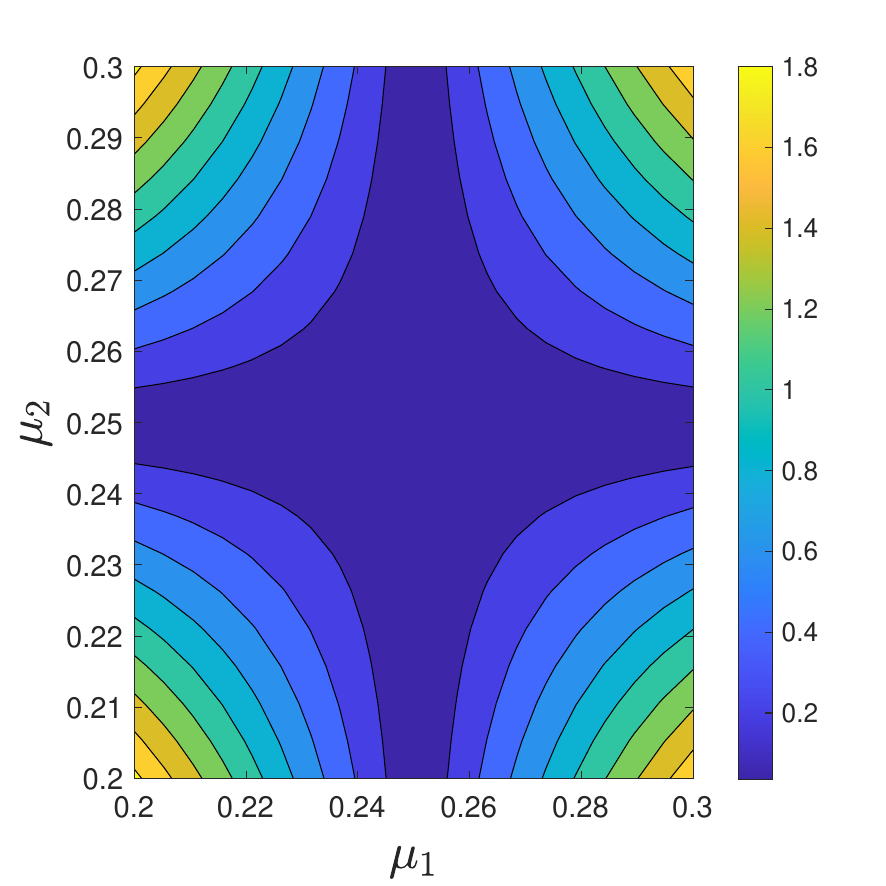}
	\label{fig10b}
	\end{subfigure}
	\caption{Percent relative error of $x^m(\mu_1,\mu_2)$ \eqref{approx} to approximate \eqref{ad-diff} for $400$ pairs of $(\mu_1,\mu_2)$, $m=50$, $n=9999$. Here accurate: $<1.5 \%$, reliable $<6 \%$. Note, same model in both figures.}
	\label{fig10}
\end{figure}

We construct approximations to \eqref{ad-diff} with Algorithm~\ref{alg:ChebyshevHOPGD}, where the sampling is performed as displayed in Figure~\ref{fig11a}. Figure~\ref{fig10} shows the percent relative error for approximations to \eqref{ad-diff} for $400$ different pairs $(\mu_1,\mu_2) \in [0,0.5]\times[0,0.5]$, where, specifically, we consider the solutions at $t=0.01$. As in \cite{HOPGD}, approximate solutions with percent relative error below $6 \%$ are considered reliable, and approximations with percent relative error below $1.5 \%$ are considered accurate. Here the error is computed by comparing the approximation to the solution obtained using backslash in \textsc{Matlab}. The $9$ snapshots used to constructed the reduced order model here were generated with $2$ executions of Chebyshev BiCG, requiring the LU decompositions of $2$ matrices of dimension $n \times n$.

Note, solutions to the advection-diffusion equation \eqref{ad-diff} are time-dependent. An all-at-once procedure could be utilized to approximate solutions to this equation at many different time-steps, though this approach would result in a larger number of unknowns and a longer simulation time. Alternatively, the decomposition could be modified to construct an approximation similar to the one in \eqref{tensor}, separated in the time variable $t$ as well as the parameters $\mu_1$ and $\mu_2$. Such an approach would, however, require specific testing. Furthermore, a decomposition of this form performed with HOPGD sampling on a sparse grid in the parameter space may not be optimal \cite{SparseHOPGD}. 

\subsection{Second simulation, parameter estimation with snapshots sampled on a sparse grid}
Consider again a parameter estimation problem. Specifically, we have an approximation to $u(x,t) \in \mathbb{R}^n$ in \eqref{ad-diff}, where $t=0.01$ is known, but the parameters $\mu_1$ and $\mu_2$ are not. This approximation corresponds to the solution to the discretized problem and is obtained using backslash in \textsc{Matlab}. Additionally, there is some uncertainty in the measurement of the given solution. In this way, we express our observed solution $u^{\text{obs}} \in \mathbb{R}^n$ as  
\begin{align} \label{obs-solution-with-error}
	u^{\text{obs}} \coloneqq u(x,t) + \varepsilon \Delta, 
\end{align}
where $\Delta \in \mathcal{N}(0,I) \in \mathbb{R}^n$ is a random vector and $\varepsilon = 10^{-2}$.  

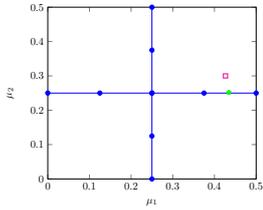
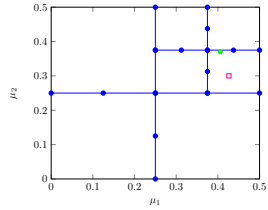
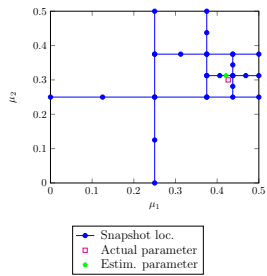
\begin{figure}[h]
\centering
		\begin{subfigure}[T]{0.45\columnwidth}
		\centering
		\begin{tikzpicture}[scale=.4]
		\begin{axis}[xlabel=$\mu_1$, ylabel=$\mu_2$,xmin=0,xmax=0.5,ymin=0,ymax=0.5,legend style={at={(.42,-.3)},
		anchor=north,legend columns=1,font=\large},legend cell align={left},each nth point={1}]

		\addplot [mark=*,color=blue]
        		table[x=column 1,y=column 2,col sep=comma]{raw-data/parest2_melina.csv};

		\addplot [mark=pentagon*,color=green]
        		table[x=column 1,y=column 2,col sep=comma]{raw-data/parest4_melina.csv};
		
		\addplot [mark=square,color=magenta]
        		table[x=column 1,y=column 2,col sep=comma]{raw-data/parest3_melina.csv};
		
		\addplot [mark=*,color=blue]
        		table[x=column 1,y=column 2,col sep=comma]{raw-data/parest1_melina.csv};

		\end{axis}
		\end{tikzpicture}
	\caption{\centering First run, $m=50$.}
	\label{fig11a}
	\end{subfigure}	
	\hspace{.4cm}
	\begin{subfigure}[T]{0.45\columnwidth}
		\centering
		\begin{tikzpicture}[scale=.4]
		\begin{axis}[xlabel=$\mu_1$, ylabel=$\mu_2$,xmin=0,xmax=0.5,ymin=0,ymax=0.5,legend style={at={(.45,-.2)},
		anchor=north,legend columns=1,font=\large},legend cell align={left},each nth point={1}]
	
		\addplot [mark=*,color=blue]
        		table[x=column 1,y=column 2,col sep=comma]{raw-data/parest6_melina.csv};
	
		\addplot [mark=*,color=blue]
        		table[x=column 1,y=column 2,col sep=comma]{raw-data/parest1_melina.csv};
	
		\addplot [mark=*,color=blue]
        		table[x=column 1,y=column 2,col sep=comma]{raw-data/parest2_melina.csv};
		
		\addplot [mark=*,color=blue]
        		table[x=column 1,y=column 2,col sep=comma]{raw-data/parest5_melina.csv};
	
		\addplot [mark=pentagon*,color=green]
        		table[x=column 1,y=column 2,col sep=comma]{raw-data/parest8_melina.csv};
		
		\addplot [mark=square,color=magenta]
        		table[x=column 1,y=column 2,col sep=comma]{raw-data/parest7_melina.csv};

		\end{axis}
		\end{tikzpicture}
	\caption{\centering Second run, $m=25$.}
	\label{fig11b}
	\end{subfigure}	
	
	\vspace{.3cm}
	
	\begin{subfigure}[T]{0.45\columnwidth}
		\centering
		\begin{tikzpicture}[scale=.4]
		\begin{axis}[xlabel=$\mu_1$, ylabel=$\mu_2$,xmin=0,xmax=0.5,ymin=0,ymax=0.5,legend style={at={(.42,-.25)},legend entries={abc},scatter/classes={
		abc={only marks}},anchor=north,legend columns=1,font=\large},legend cell align={left},each nth point={1}]
	
		\addplot [mark=*,color=blue]
        		table[x=column 1,y=column 2,col sep=comma]{raw-data/parest9_melina.csv};
		\addlegendentry{\hphantom{ }Snapshot loc.}	
		
		\addplot [mark=square,color=magenta,only marks]
        		table[x=column 1,y=column 2,col sep=comma]{raw-data/parest11_melina.csv};
		\addlegendentry{\hphantom{ }Actual parameter}
		
		\addplot [mark=pentagon*,color=green,only marks]
        		table[x=column 1,y=column 2,col sep=comma]{raw-data/parest12_melina.csv};
		\addlegendentry{\hphantom{ }Estim. parameter}
			
		\addplot [mark=*,color=blue]
        		table[x=column 1,y=column 2,col sep=comma]{raw-data/parest1_melina.csv};
	
		\addplot [mark=*,color=blue]
        		table[x=column 1,y=column 2,col sep=comma]{raw-data/parest2_melina.csv};
		
		\addplot [mark=*,color=blue]
        		table[x=column 1,y=column 2,col sep=comma]{raw-data/parest5_melina.csv};
	
		\addplot [mark=*,color=blue]
        		table[x=column 1,y=column 2,col sep=comma]{raw-data/parest6_melina.csv};
		
		\addplot [mark=*,color=blue]
        		table[x=column 1,y=column 2,col sep=comma]{raw-data/parest10_melina.csv};
	
		\end{axis}
		\end{tikzpicture}
	\caption{\centering Third run, $m=15$.}
	\label{fig11c}
	\end{subfigure}	
\caption{Parameter estimation for noisy solution to \eqref{ad-diff}, $n=9999$.}
\label{fig11}
\end{figure}

Figure~\ref{fig11} and Table~\ref{table:2} show the result of this parameter estimation problem, using a similar strategy to the one described in \cite{HOPGD3}. Here $3$ successive HOPGD models are constructed from snapshots sampled on a sparse grid in the parameter space $(\mu_1,\mu_2) \in [0,0.5]\times[0,0.5]$. This simulation executes a modified version of Chebyshev BiCG $6$ times, generating the $23$ snapshots efficiently. Note, this requires the LU decompositions of $6$ matrices of dimension $n \times n$. After the third run of Algorithm~\ref{alg:ChebyshevHOPGD}, the estimated parameter is very close to the actual parameter. Note, the relative error in $\mu_1$ and $\mu_2$ after the third run of the simulation is of the same order of magnitude as the noise in the observed solution \eqref{obs-solution-with-error}. 

\begin{table}[H]
\centering
\begin{center}
\begin{tabular}{ c | r | r }
 & rel err $\%$ $\mu_1$ & rel err $\%$ $\mu_2$ \\ \hline
First run, $m=50$ & 1.9004 & 16.1248 \\ \hline  
Second run, $m=25$ & 4.8542 & 23.2504 \\ \hline  
Third run, $m=15$ & 1.2313 & 4.2588 \\ \hline
\end{tabular}
\end{center}
\caption{Percent relative error for parameter estimation in $\mu_1$ and $\mu_2$ for noisy solution to~\eqref{ad-diff}, $n=9999$.}
\label{table:2}
\end{table}

\section{Conclusions and future work}
This work overcomes the challenge of treating linear systems with dependence on multiple parameters efficiently.  By
generating a reduced order model for the solution using snapshots generated by carefully varying a single parameter at a time, we
are able to propose a method built on the robust method for treating systems with a single-parameter dependence
\cite{CorrentyEtAl2}, as we can use preconditioned Chebyshev BiCG to generate many snapshots simultaneously. This allows
us to to generate the snapshot solutions required to build a reduced order model in a way that is both low-cost and novel.
The model is based on an efficient decomposition of a tensor matrix, where the sampling is performed on a sparse grid.  

Tensor decompositions may fail to reach a certain level of accuracy on a given set of snapshot solutions. Further, reduced order models based on sampling can suffer from overfitting and underfitting, and the choice of snapshots can affect the simulation time. Our approach offers a way to generate a new set of snapshots on the same parameter space with little extra computation. This is advantageous, as it is not known a priori if a given set of snapshots will produce an accurate and robust model, and generating snapshots is computationally demanding in general. The reduced order model can also be used to solve a parameter estimation problem. Numerical simulations show competitive results.  

In \cite{SparseHOPGD}, a residual-based accelerator was used in order to decrease the number of iterations required by the fixed point method when computing the updates described in \eqref{phi-m}, \eqref{sol-f}, and \eqref{sol-f2}. Techniques of this variety are especially effective when the snapshot solutions have a strong dependence on the parameters. Such a strategy could be directly incorporated into this work, though specific testing of the effectiveness would be required. 

\section{Acknowledgements}
The Erasmus\raisebox{0.25ex}{\scriptsize +} programme of the European Union funded the first author's extended visit to Trinity College Dublin to carry out this research. The authors thank Elias Jarlebring (KTH Royal Institute of Technology) for fruitful discussions and for providing feedback on the manuscript. Additionally, the authors wish to thank Anna-Karin Tornberg (KTH Royal Institute of Technology) and her research group for many supportive, constructive discussions.

The authors declared that they have no conflict of interest.

\bibliographystyle{siam}
\bibliography{siobhanbib,eliasbib}

\end{document}